\documentclass[numreferences]{kluwer}

\usepackage{amsfonts,amsthm}
\usepackage{psfig,graphicx,psfrag}
\usepackage{klucite}

\usepackage{amssymb}

\def\SS{\mathbb{S}}

\def\IH{{{\mathbb{H}}^3}}

\def\HHH{\mathcal{H}}
\def\cal{\mathcal}
\def\RR{\mathbb{R}}

\def\D{\Delta}

\def\Di{\Delta_{int}}
\def\Ni{N_{int}}

\def\De{\Delta_{ext}}
\def\Ne{N_{ext}}

\def\PP{\mathbb{P}}
\def\MM{\mathbb{M}}

\def\oC{{\hat{\mathbb {C}}}}

\def\CCC{\mathcal{C}}

\def\IC{\Bbb C}

\def\De{\Delta_{ext}}
\def\D{\Delta}

\def\C{{\rm{ Circ }} }

\def\x(AB)i{{Ax_{{A \cdot B^{-1}}}}}
\def\(AB)i{{{{A \cdot B^{-1}}}}}

\def\T{{Tr \;}}
\def\bB{{\bar{B}}}

\def\pl{pl}

\def\PPi{{\mathbb{P}}_{int}}
\def\PPe{{\mathbb{P}}_{ext}}
\newtheorem{thm}{Theorem}[section]

\newtheorem{lemma}[thm]{Lemma}

\newtheorem{cor}[thm]{Corollary}

\newtheorem{prop}[thm]{Proposition}

\begin{document}
\begin{article}
\begin{opening}
\title{The Geometry of Two Generator Groups:\\Hyperelliptic Handlebodies}

\author{Jane \surname{Gilman}\email{gilman@andromeda.rutgers.edu}\thanks{Supported in part by NSA
grant \#MSPF-02G-186 and a grant from the Rutgers Research Council
and thanks Kings College, University of London for its hospitality
during the preparation of some of
 this work.}}\institute{ Department of Mathematics, Smith Hall\\
 Rutgers University, Newark, NJ 07102, USA \\ }
\author{Linda \surname{Keen}\email{linda.keen@lehman.cuny.edu}\thanks{Supported in
part by a PSC-CUNY grant.}} \institute{ Mathematics
Department\\CUNY Lehman College and Graduate Center\\ Bronx, NY 10468, USA}

\runningtitle{Hyperelliptic Handlebodies} \runningauthor{J. Gilman
and L. Keen}
\begin{abstract} A Kleinian group  naturally
stabilizes certain subdomains and closed subsets of the closure of
hyperbolic three space and yields a number of different quotient
surfaces and manifolds. Some of these quotients have conformal
structures and others hyperbolic structures. For two generator
free Fuchsian
 groups,
 the quotient three manifold is
a genus two solid handlebody
 and its boundary is a hyperelliptic Riemann surface. The convex
 core is also a hyperelliptic Riemann surface.  We find
 the Weierstrass points of both of these surfaces.
We  then generalize the notion of a hyperelliptic Riemann surface
to a ``hyperelliptic'' three manifold. We show that the handlebody
has a unique order two isometry fixing six unique geodesic line
segments, which we call the {\sl Weierstrass lines} of the
handlebody. The Weierstrass lines are, of course, the analogue of
the Weierstrass points on the boundary surface.
 Further, we show that the manifold is
foliated by surfaces equidistant from the convex core, each fixed
by the isometry of order two. The restriction of this involution
to the equidistant surface fixes six {\sl generalized Weierstrass
points} on the surface. In addition, on each of these
 equidistant surfaces we
 find an orientation reversing involution that
fixes curves through the {\sl generalized Weierstrass points}.
\end{abstract}




\end{opening}


\section{Introduction}
We begin with a Kleinian group, $G$, a discrete group of M\"obius
transformations. We can think of this group as acting on the
closure of hyperbolic three space. That is, on the union of $\IH$
hyperbolic three space and its boundary, $\oC$ the complex sphere.
In addition to stabilizing each of these, the group naturally
stabilizes certain subdomains and closed subsets, yielding  a
number of different quotient surfaces and manifolds. Some of these
quotients have conformal structures and others hyperbolic
structures.

If the  group is a Fuchsian group,  it fixes both the interior and
the exterior of a  disc in the complex plane, and  a corresponding
 hyperbolic plane lying over the
 disc.  The limit set of the group is contained in the circle separating
 the  discs in the plane.
 In this  case there are
additional quotient surfaces to consider and certain symmetries
within the  quotients or among the quotients.

Kleinian groups whose limit sets are proper subsets of $\oC$ and
Fuchsian groups whose limit sets are proper subsets of the
boundary of their invariant disc are known as ``groups of the
second kind''. In this paper, we consider only finitely generated
 groups of the second kind that are also free groups.
 For such
 groups,
we find the relations between all of these quotients and the
relations between the various hyperbolic and conformal metrics
that can be placed on these spaces.

The quotient surfaces include the Nielsen kernel, the Fuchsian
 quotient and its  Schottky double. There is also a
 Nielsen double. The Schottky double is the boundary of the  quotient three
manifold. For a two generator Fuchsian group (of the second kind)
we prove that the boundary of the convex core of the three
manifold is the image under a pleating map of the Nielsen  double.
That is, we find the Fuchsian group that uniformizes the  Nielsen
double and use it to obtain an isometric map from the Nielsen
double to the boundary of the convex core of the manifold. This
map is a pleating map.

We are interested primarily in the case where the groups are
two-generator groups. For these groups (when they are free and
discrete) the quotient three manifold is a genus two solid
handlebody
 and its boundary is a hyperelliptic Riemann surface.
In this case we  generalize the notion of a hyperelliptic Riemann
surface to a ``hyperelliptic'' three manifold. We show that the
handlebody has a unique order two isometry fixing six unique
geodesic line segments, which we call the {\sl Weierstrass lines}
of the handlebody. The Weierstrass lines are, of course, the
analogue of the Weierstrass points on the boundary surface.
 Further, we show that the manifold is
foliated by surfaces equidistant from the convex core, each fixed
by the isometry of order two. The restriction of this involution
to the equidistant surface fixes six {\sl generalized Weierstrass
points} on the surface. In addition, on each of these
 equidistant surfaces we
 find an orientation reversing involution that
fixes curves through the generalized Weierstrass points.
\vskip .1in The organization of the paper is outlined in the table
of contents below. We begin with  basic definitions and  notation.
We then define the various quotient spaces of the Kleinian group;
that is, the Fuchsian surfaces, the Schottky and Nielsen doubles,
the handlebody and the convex core. It turns out that there is a
natural dichotomy based on whether the axes of the generators are
either disjoint or intersect.  We carry out our discussion in two
parallel parts; in part 1, we state and prove all our results in
the disjoint axes case and in part 2 we do the same for the
intersecting axes case.

 \tableofcontents
 \nopagebreak
\section{Setup and Notation}
\label{sec:notation}
 We let $G$ be a finitely generated discrete
group of M\"obius transformations.
\subsection{M\"obius  transformations and metrics}

A fractional linear transformation or a M{\"o}bius transformation
is a conformal homeomorphism of $\oC$ of the form $z \rightarrow
{\frac{az+b}{cz+d}}$ where $a,b,c,d \in \IC,$  $ad-bc = 1$. We let
$\MM$ be the group of M{\"o}bius transformations. We also consider
elements of $\overline{\MM}$, the group of anti-conformal
M{\"o}bius transformations. These act on $\oC$ as fractional
linear transformations sending $z \rightarrow
{\frac{a{\bar{z}}+b}{c{\bar{z}}+d}}$.

We use the upper-half-space model for ${\IH}$. If the coordinates
for Euclidean three-space are $(x,y,t)$, then $\IH =\{ (x,y,t)\; |\;
t > 0 \}$; the plane with $t=0$ is $\partial \IH = \oC$.

A Euclidean sphere whose center is on $\oC$, intersects $\oC$ in a
circle and $\IH$ in  Euclidean hemisphere whose {\sl horizon} is
that circle. Inversion in such a Euclidean sphere fixes  $\IH$. It
also fixes $\oC$ and its restriction to $\oC$ has the same action
as inversion in the horizon. We include Euclidean planes
perpendicular to $\oC$ as spheres, thinking of them as ``spheres''
through $\infty$; the horizon of such a plane is a straight line
or ``circle'' through $\infty$ on $\oC$. With this convention the
terms {\sl inversion} and {\sl reflection} are used
interchangeably.

The action of elements of $\MM$ and $\overline{\MM}$ on $\IH$ is
related to the action on $\oC$ as follows:  each element of $\MM$
can be factored as a product of reflections in circles lying in
$\oC$ and its action on $\oC$ is the restriction to $\oC$ of the
product of reflections in the Euclidean hemispheres whose horizons
are these circles.

There is a natural metric on $\IH$, the hyperbolic metric,
preserved by  these reflections.  Products of even numbers of
reflections are the orientation preserving isometries for this
metric and the hyperbolic metric is the unique metric for which
$\MM$ is the full group of orientation preserving isometries.
Geodesic lines in this metric are circles (and straight lines)
orthogonal to $\oC$. The Euclidean hemispheres are geodesic
surfaces called {\em hyperbolic planes}. The boundary of a
hyperbolic plane in $\oC$ is called the horizon of the plane since
it is at infinite distance from every point in the plane.

We set
$$\Di= \{z: |z|<1\}, \;\;\De=\{z: |z|>1\}\;\; \mbox{  and } \;\;
\partial\Delta = \{z: |z|=1\}.$$ For readability we use $\D$ for
$\Di$ when no confusion is apt to arise.

There is also a natural hyperbolic metric on $\D$ viewed as the
hyperbolic plane. The isometries are elements of $\MM$ that fix
$\D$ and geodesics are circles orthogonal to $\partial\D$.

We let $\PP$ denote the plane in $\IH$ whose horizon is the unit
circle. We can think of $\PP$ as consisting of $\PPe$ and $\PPi$,
the hyperplane boundaries of hyperbolic half-spaces  lying over
$\De$ and $\Di$ respectively. The hyperbolic metric in $\IH$
 restricts to a metric on $\PP$ and
orthogonal projection from $\PPi$ (or $\PPe$) with this metric to
$\Di$ (or $\De$) with the usual hyperbolic metric $\rho_{\D}$ (or
$\rho_{\De}$) is an isometry.

We use the notation $R_{\partial\Delta}$ to denote the reflection
 in $\partial \Delta$ mapping $\oC$ to itself and $R_{\PP}$ to denote
the reflection in $\PP$ mapping $\IH$ to itself.

We remind the reader that elements of $\MM$ have both an algebraic
 classification by their traces, and corresponding geometric classification
according to their action on $\IH$ or $\oC$.

In this paper, we are  concerned with finitely generated free
groups that are subgroups of $\MM$. They contain no elliptic
elements. We shall also make the simplifying assumption that there
are no parabolic elements.   Every element $A \in \MM$ therefore
has two fixed points in $\oC$.  The line $Ax_A$ in $\IH$ joining
the fixed points is invariant under $A$ and is called its axis.
Similarly, if the fixed points of $A$ lie on $\partial\D$, the
geodesic in $\D$ joining them is also called the axis.  We will
also denote this axis by $Ax_A$.  The context will make clear which
axis we mean.

We will also have occasion to consider certain groups generated by
elliptic elements of order $2$. These groups contain free groups
as subgroups of index $2$ as well as subgroups containing  certain
orientation reversing elements of order two.

\subsection{Hyperelliptic surfaces and the hyperlliptic involution}
\label{sec:involutions}

If a compact Riemann surface of genus $g$  admits a conformal
involution with $2g+2$ fixed points, the involution is unique and
the surface  is called a {\sl hyperelliptic surface}. The fixed
points of the hyperelliptic involution are called the {\sl
Weierstrass points} (\cite{Spr} p. 275 and \cite{Hur}) of the
surface. In addition to this geometric definition, these points
also have  an analytic definition in terms of the possible zeros
and poles of meromorphic functions on the surface. Every Riemann
surface of genus two is hyperelliptic. It  admits a unique
conformal involution with six  fixed points,  the {\sl
hyperelliptic involution}.

\section{Subspaces and Quotient Spaces}
\label{sec:quotient spaces}
\subsection{Limit sets, regular sets, convex regions and their quotients}
\label{sec:limitsetsetc}
 We assume that $G$ is a finitely generated
Kleinian group. Since it is discrete, it acts discontinuously
everywhere in $\IH$. We denote the subset of $\oC$ where it acts
discontinuously by
 $\Omega(G)$ and denote its complement in $\oC$, which is called the {\em limit set}, by
  $\Lambda(G)=\oC \setminus \Omega(G) $.

 A {\em fundamental domain} $F \subset \SS$ for $G$ acting discontinuously on a space
 $\SS=\Delta,\PP,\Omega, \ldots$ is an open set such that
 if $z_1,z_2$ are a pair of points in the interior of $F$ and if
$W(z_1)=W(z_2)$ for some $W \in G$, then $W=Id$ and $\cup_{W \in
G}\overline{W(F)}=\SS$.

If $W$ is any closed subset of $\oC$, we call the smallest closed
hyperbolically convex subset of $\IH$ containing it, the {\sl
convex hull of $W$ in $\IH$}. If $W$ is the limit set of $G$ we
denote its convex hull by $\CCC(G)$ or simply $\CCC$. Analogously, if $W$ is
any closed subset of $\partial{\D}$, we call the smallest closed
hyperbolically convex subset of $\D$ containing it the {\sl convex hull of $W$
in $\D$}.

Since the set $W$ we use in this paper is  the limit set $\Lambda$
of a Fuchsian group, we need to distinguish between its convex
hull in $\IH$ and its convex hull in $\D$.  Because  the convex
hull in $\D$ of the limit set of a Fuchsian group historically was
called the {\sl Nielsen (convex) region for $G$}, we continue that
tradition and denote it by $K(G)$ or $K$.

\vspace{.2in}
 \centerline{\bf The Quotients for free Kleinian groups  of the second kind }

\begin{description}
\item[The surface $S$]: As $G$ is a group of the second kind,
$\Omega$ is connected.  The holomorphic projection $\Omega(G)
\rightarrow \Omega(G)/G=S$ determines  a Riemann surface $S$.
Since $G$ is finitely generated, free, and contains no parabolics,
$S$ is a compact surface of finite genus.
\end{description}

\begin{description}
  \item[The solid handlebody $\bar{H}$]:  The projection $\IH \rightarrow
\IH/G=H$ determines $H$ as a complete hyperbolic 3-manifold.  Its
boundary is the Riemann surface $S = \Omega(G)/G$. That is,
$\partial H=S$. Therefore $\bar{H}$ is a solid handlebody.
\end{description}

\begin{description}
\item[The Convex Core $\CCC/G$]: The projection $\CCC(G) \rightarrow
\CCC(G)/G=\CCC/G$ is a convex closed submanifold of $H$. Its
boundary, $\partial\CCC/G$ is homeomorphic to the surface $S$.
\end{description}

\vspace{.1in}
 \centerline{\bf The Quotients for free Fuchsian groups  of the second kind }

 \vspace{.1in}
If $G$ is a Fuchsian group, replacing $G$ by a conjugate if
necessary we may assume that $G$ acts invariantly on the unit
disc.
\begin{description}
\item[The Fuchsian surface $\D/G$]: The projection $\D \rightarrow
\D/G$ is a complete Riemann surface that we call the {\em
Fuchsian surface of $G$}. Since $G$ is free and contains no
parabolics, $\D/G$ is a surface of finite genus with finitely many
removed discs.
\end{description}


\begin{description}
\item[The Nielsen kernel $K(G)/G$]: The {\sl Nielsen kernel of $G$}
 is  the quotient of the Nielsen region;
 that is, the surface $K(G)/G$.
  It is also called the Nielsen kernel of the surface $\D/G$.  It is homeomorphic to $\D/G$.
\end{description}

\vspace{.1in}

If $G$ is a free group of rank two, then as we will see below, $S$
is a compact surface of genus two. Thus it  is the boundary of  a
handlebody $H$ of genus two.

\subsection{Schottky doubles, funnels and Nielsen doubles}
\label{sec:doubles}
 In this section we use the  exposition of \cite{Bers}
for  constructing  the Schottky double and the Nielsen double for
Fuchsian groups of the second kind. We also define metrics for
them.

For any compact Riemann surface of genus $g$ with $k$ holes, one
can always form the conformal double.  We denote the conformal
double of $\D/G$ by $S$ and the conformal double of $K(G)/G$ by
$S_K$. Bers termed $S$  the {\sl Schottky double} and clarified
the relationship  between $S$ and $S_K$. (See
theorem~\ref{thm:bers} below.) To distinguish between these
doubles we call $S_K$ the {\em Nielsen double}.

The anti-conformal reflection in $\partial \D$ induces a
canonical anti-conformal involution on $S$ which we denote by $J$.

As a compact Riemann surface, $S$ admits a uniformization by a
Fuchsian group $\hat{G}$ and the projection $\D \rightarrow
\D/\hat{G}=S$ defines a hyperbolic metric on $S$.  The restriction
of this metric to the Fuchsian surface $\D/G \subset S$ is called
the {\em intrinsic metric} on $\D/G$.

The Nielsen kernel of $S$, $K/G$, is a Riemann surface
homeomeorphic to $S$. One can construct its conformal double $S_K$
which we call the {\em Nielsen double}. Again the restriction of
the hyperbolic structure on $S_K$ defines the intrinsic metric on
the Nielsen kernel. The intrinsic metric agrees with the
restriction of the hyperbolic metric on the Fuchsian surface
$\D/G$.

The intersection of $\Omega(G)$ with $\partial \Delta$ consists of
an infinite set $\beta$ of open intervals $I$, called the {\sl
intervals of discontinuity of $G$}. The images of the intervals of
discontinuity project to (ideal) boundary curves of $\D/G$.

Each $I \in \beta$ corresponds to the axis of some element in $A_I
\in G$. Let $\C_{A_I}$ be the circle determined by the axis of
$A_I$ (that is, the circle in $\oC$ orthogonal to $\partial\Delta$
through the fixed points of $A_I$), and let $D(A_I)$ be the disc
it  bounds containing $I$. The quotient  $F(A_I) = D(A_I)/G
\subset S$ is called a {\em funnel} in $S$.   The curves
 $\C_{A_I} \cap \De$,  $\C_{A_I} \cap \Di$ and $I$
project  onto homotopic curves on $S$. The first two are the {\em
boundary curves of the funnel} and the last is the {\em central
curve of the funnel}.

The following theorem follows from results in \cite{Bers}
\begin{thm}\label{thm:bers} We have the following relationship
between the Nielsen kernel and the Schottky double:

$$S=\Omega(G)/G = K(G)/G \cup J(K(G)/G) \cup_{I \in \beta} F(A_I)$$

Moreover, the central curves of the funnels are geodesics on  $S$.
The involution $J$ interchanges the boundary curves and fixes the
central curve.
\end{thm}

\subsection{The three-manifold}
\label{sec:three manifold}
 We continue assuming $G$ a Fuchsian
group of the second kind with invariant disc $\D$. As a group
acting on $\IH$ it acts invariantly on a hyperbolic plane $\PP$
whose horizon is $\partial \D$. In this section, we consider it as
acting on $\PP$.

\begin{description}
\item[The  Nielsen double, ${S}_{K_{\PP}}$]: The Nielsen region for
$G$ acting on $\PP$, $K_{\PP}$, is the convex hull in $\PP$ of
$\Lambda(G)$. If we think of $\PP$ as having two sides, $\PPi$ and
$\PPe$ we form the doubled kernel by taking one copy on each side:
${S}_{K_{\PP}} = K_{\PPi}/G \cup K_{\PPe}/G$. This is
isometrically equivalent to the  Nielsen double $S_K$ defined in
the previous section.

 \item[The convex hull, $\CCC$]: The hyperbolic convex hull of $\Lambda(G)$ in $\IH$  is
$\CCC=\CCC(G)$. Because $G$ is Fuchsian, $\CCC$ is contained in
$\PP$ and has no interior so that
 $\partial\CCC=\CCC$.  It is, therefore, sometimes natural to think of
$\partial\CCC$ as two-sided with $\CCC_{int}$, the interior side
facing $\Di$ and $\CCC_{ext}$, the exterior side facing $\De$.
Note we have the obvious identifications $\CCC=K_{\PP}$,
$\CCC_{int}=K_{\PPi}$ and $\CCC_{ext}=K_{\PPe}$.
 \item[The convex core of $H$, $N= \Ni \cup \Ne$]: Since $\CCC=\partial\CCC$  the convex core
$N=\CCC/G$ can be thought of as consisting of two hyperbolic
surfaces $N_{int}=\CCC_{int}/G$ and $N_{ext}=\CCC_{ext}/G$ joined
along their common boundary curves. Again we have the obvious
identifications $\Ni=K_{\PPi}/G$, $\Ne=K_{\PPe}/G$. Note that $N$
is isometric to ${S}_{K_{\PP}}$.
\end{description}
We conclude this section with a concept that relates surfaces and
three manifolds.

\begin{description}
\item [The Pleated surface $({\tilde{S}}, pl)$]: A pair,
$({\tilde{S}}, pl)$, is a {\sl  pleated surface}, if $\tilde{S}$
is a complete hyperbolic surface and $\pl$ is an immersion and  a
hyperbolic isometry of $\tilde{S}$ into a hyperbolic three
manifold $H$, such that every point in $\tilde{S}$ is in the
interior of some geodesic arc which is mapped to a geodesic in
$H$. The {\sl pleating locus} is  the set of points in $\tilde{S}$
that are contained in the interior of exactly one geodesic arc
that is mapped by $pl$ to a geodesic arc.
\end{description}

Our  goal is to explain the connection between the various
quotient Riemann surfaces and related quotient hyperbolic surfaces
and quotient three manifolds we have now defined.  We will
restrict our discussion to the situation where
  $G$ is a free two generator
Fuchsian group.

We will show that if $G$ is a two generator Fuchsian group, then
$N$, with its intrinsic metric, is the image of the  Nielsen
double $S_K$, with its intrinsic metric, under a pleating map,
$pl$. We will see that the image of the pleating locus  consists
of the identified
 geodesic boundary curves of $N_{int}$ and $N_{ext}$.
Further we will show that $\bar{H} = {\cal{S}\times I}$ where
$\cal{S}$ is a surface of genus $2$, $I=[0,\infty]$ and
${\cal{S}\times 0}=N$, ${\cal{S}\times \infty}=S$, and the
surfaces  ${\cal{S}\times s}=S(s)$ are {\sl equidistant surfaces}
from the convex core $N$.

We begin with the discrete free Fuchsian group on two generators,
 $G= <C,D>$.  We want to find a
 set of generators $\langle A,B \rangle$ for $G$ that have two special properties:
 first, they are {\em geometric},
 that is, they are the side pairings for a fundamental polygon (\cite{Pu}) and second, the
 traces  of $A$ and $B$ are less than the traces for any other pair of generators.
 Although we know $G$ is discrete, we apply what is known
as the discreteness algorithm to find the special generators. They
are the pair of generators at which the algorithm stops. We thus
 call  $\langle A,B \rangle$ the {\em stopping generators} for $G$.

 The first step of the algorithm determines whether the axes of the
generators intersect
 or not.
This is determined by computing the trace of the commutator,
$[C,D]=CDC^{-1}D^{-1}$  (see \cite{Gint} or \cite{RR}).
 If $\T [C,D] <-2,$ the axes intersect; if $ -2 < \T [C,D]
< 2,$ the group contains an elliptic element; and if $\T [C,D]
>2,$ the axes are disjoint.  Either the axes intersect for every
Nielsen equivalent pair of generators or they are disjoint for
every pair.  Since $G$ is free, it contains no elliptics and since
we have assumed it contains no parabolics, we have $|\T [C,D]|>2$.

Our constructions are somewhat different depending on whether the
axes of the generators are disjoint or not.  We separate these
constructions into two parts.

\nopagebreak
\part{DISJOINT AXES: $\T [C,D] > 2$} \nopagebreak
\label{part:part1} \vspace{0.1in}

Throughout part I we assume the axes of the generators of $G$ are
disjoint. We begin by finding fundamental domains for $G$ acting
on $\D$ and on $K(G)$.  Next, we find the Weierstrass points for
$S$.  We then turn to $\IH$ and find fundamental domains for $G$
acting on $\PP$ and on $\CCC$. We construct the Fuchsian group
uniformizing $S_K$ and the pleating map sending $K(G)$ onto
$\partial{\CCC}$.  We then use the pleating map to find the
Weierstrass points of $S_K$.  Finally, we work with the handlebody
$H$ and construct the Weierstrass lines and the equidistant
surfaces.  We conclude with a discussion of  anti-conformal
involutions on the various quotients.

\section{Subspaces of $\Delta$ and quotients by $G$}

\label{sec:schottkynonint} \subsection{Fundamental Domains in $\D$
and $K(G)$} \label{sec:funddomains non-int}

 When the axes of the generators of $G$ are disjoint, the Gilman-Maskit
  discreteness algorithm \cite{GilMas} proceeds by successively
  replacing the generators with a Nielsen equivalent pair.  After
  finitely
 many steps, it decides whether the group generated by the original generators is discrete or not.
   If the group is discrete,
 the final set of generators,  the stopping generators, are both
 geometric and have minimal traces.

 The reflection lines
$(L,L_A,L_B)$ defined as follows for any pair of generators, play
a major role in the algorithm:
\begin{quote} $L$ is the common perpendicular in $\Delta$ to the axes of $A$
and $B$\\
$L_A$ is the geodesic in $\Delta$ such that $A$ factors as the
product $A=R_LR_{L_A}$ where $R_K$ denotes reflection in the
geodesic $K$. \\
$L_B$ is the geodesic in $\Delta$ such that $B=R_LR_{L_B}$
\end{quote}

The stopping condition for the algorithm is that the lines
$(L,L_A,L_B)$ bound a domain in $\Delta$; that is, no one
separates the other two. This is equivalent to $R_L, R_{L_A},
R_{L_B}$ being geometric generators for the index two extension of
$G$, $\langle G, R_L \rangle$;  the domain bounded by
$(L,L_A,L_B)$ is its fundamental domain.

Note that $A^{-1}B=R_{L_A}R_{L_B}$ so that its axis is the common
perpendicular to $L_A$ and $L_B$.  In fact, the axes $(Ax_A, Ax_B,
Ax_{A^{-1}B})$ bound a domain if and only if the lines
$(L,L_A,L_B)$ do.  We orient the lines so that $L$ is
 oriented from $Ax_A$ to $Ax_B$, $L_B$ is oriented from
 $Ax_B$ to $Ax_{A^{-1}B}$ and $L_A$ is oriented from
 $Ax_{A^{-1}B}$ to $Ax_A$.

Taking the union of the domain bounded by
  the lines $(L,L_A,L_B)$ together with its reflection in $L$, we
 obtain a domain $F$ in $\D$ bounded by the lines  $L_A,$
 $L_B$, $L_{\bar{A}}=R_L(L_A)=A(L_A)$ and
 $L_{\bar{B}}=R_L(L_B)=B(L_B)$.

\begin{figure}[hbt]\centering

\psfrag{LbB}[cc]{$L_{\bar{B}}$}
 \psfrag{LbA}[cc]{$L_{\bar{A}}$}
\psfrag{PL}[cc]{$p_L$}
 \psfrag{L}[cc]{$L$}
\psfrag{QLA}[cc]{$q_{L_A}$} \psfrag{LA}[cc]{$L_A$}
\psfrag{PLA}[cc]{$p_{L_A}$}
 \psfrag{QLB}[cc]{$q_{L_B}$}
\psfrag{LB}[cc]{$L_{B}$}
 \psfrag{QL}[cc]{$q_L$}
\psfrag{PLB}[cc]{$p_{L_B}$}
 \psfrag{Delta}[cc]{$\Delta$}

\includegraphics[width=5in]{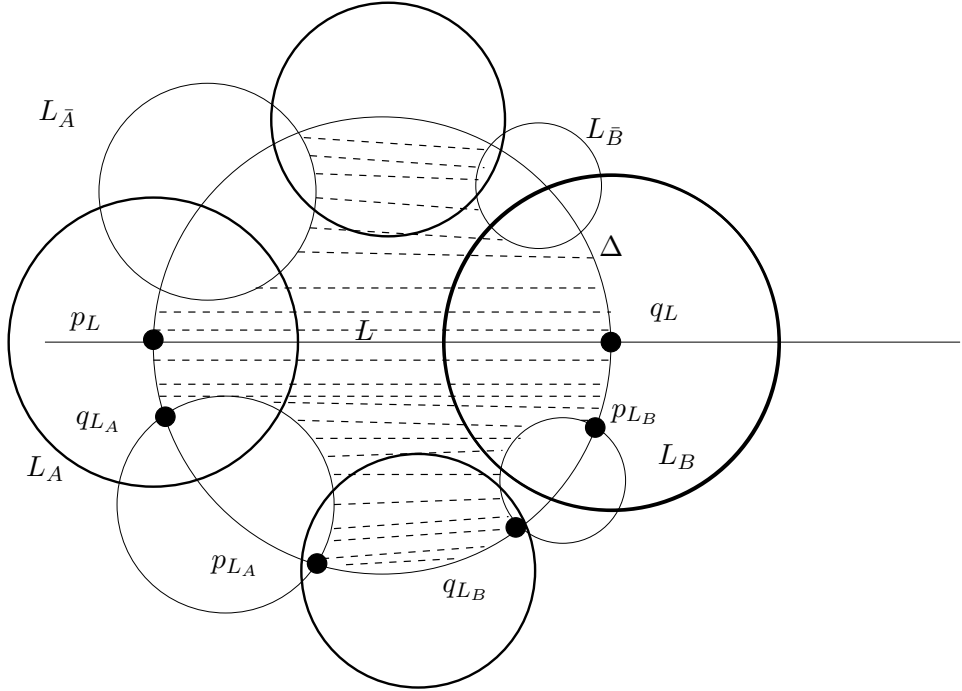}

\caption{The domain $F \subset \D$ is shaded. The axes, drawn with
their reflections in $\partial\D$ are the darker circles; the
lines $L,L_A$, $L_B,L_{\bar{A}},L_{\bar B}$ drawn with their
reflections are the lighter circles. \label{fig:s7}}
\end{figure}

\begin{figure}[hbt]\centering
\psfrag{LbB}[cc]{$L_{\bar{B}}$}
 \psfrag{LbA}[cc]{$L_{\bar{A}}$}
\psfrag{PL}[cc]{$p_L$} \psfrag{L}[cc]{$L$}
 \psfrag{LA}[cc]{$\;$ \linebreak
$\;$\linebreak $\;$\linebreak $\;$\linebreak $\;$\linebreak
$\;$\linebreak $\;$\linebreak $\;$\linebreak $\;$ \linebreak
$\;\;\;\;\;\;\;\;\;\;$  \linebreak $L_A$}
\psfrag{PLA}[cc]{$p_{L_A}$} \psfrag{QLB}[cc]{$q_{L_B}$}
\psfrag{LB}[cc]{$L_{B}$} \psfrag{QL}[cc]{$L_B$}
\psfrag{PLB}[cc]{$p_{L_B}$}
 \psfrag{Delta}[cc]{$\Delta$}
\includegraphics[width=5in]{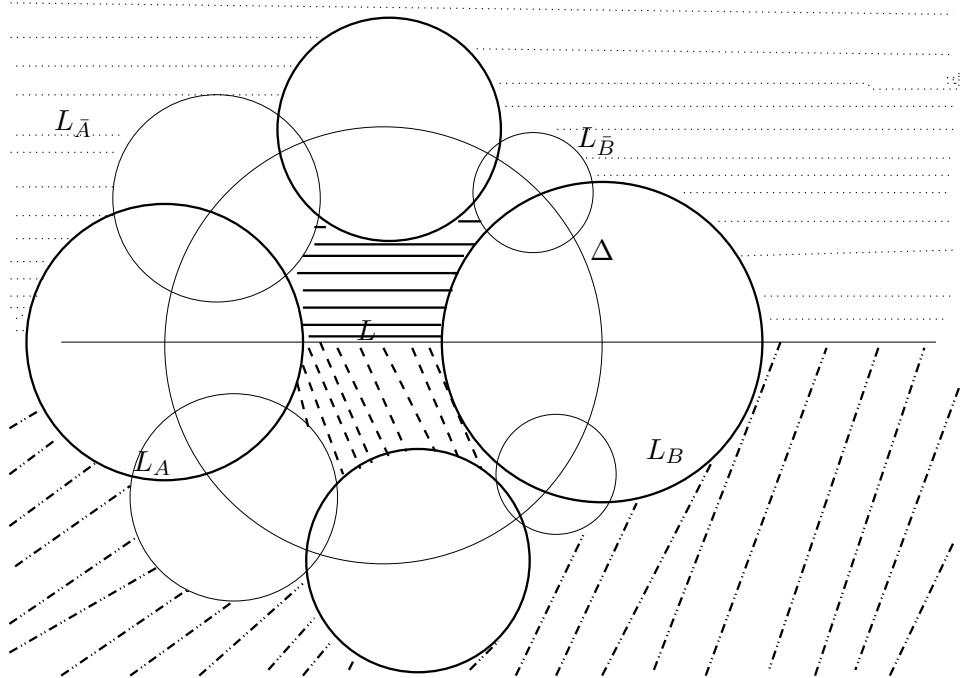}
\caption{The domain $F_K \subset \D$ and its reflection
$R_{\partial\D}(F_K)$ are shaded \label{figure:s6}}
\end{figure}

%
%
%
%
%
%
%
%
%
%
%
%
%
%
%

 \begin{prop}\label{prop:funddomain nonint} If $\langle A,B \rangle$ are
 the stopping generators for $G$, then $F$ is a fundamental domain
 for $G$ acting in $\D$. If $F_{K}$ is the domain obtained from
 $F$ by truncating along the axes of $A,B,A^{-1}B$ and $AB^{-1}$
 then $F_K$ is a fundamental domain for the action of $G$ on the
 Nielsen region of $G$ in $\D$.
 \end{prop}
 \begin{proof} (See figures 1 and 2) By construction, the algorithm stops at the pair
 of generators $\langle A,B \rangle$ such that $F$ is a
 fundamental domain for $G$ acting on $\D$ (\cite{GilMas}).

 Let $D(I_A)$ be the half plane bounded by $Ax_A$ that does
 not contain the segment of $L$ between the axes of $A$ and $B$.
 Note that it is invariant under $A$.
  Tiling $\D$ with copies of $F$, we
see that the only images of $F$ that intersect $D(I_A)$ are of the
form $A^n(F)$, for some integer $n$.  It follows that $D(I_A)$ is
stabilized by the cyclic group $\langle A \rangle$ and the
interval $I_A$ joining the fixed points of $Ax_{A}$ is an interval
of discontinuity for $G$.  The same is true for the other three
intervals $I_B,I_{A^{-1}B}$ and $I_{AB^{-1}}$ corresponding to the
other three  axes that intersect $F$.   We conclude that these
axes all lie on the boundary of the Nielsen convex region for $G$.

  Identifying the sides of $F$ we easily see that
  $\D/G$ is a sphere with three holes. The three ideal boundary
  curves are the projections of the intervals $I_A,I_{B}$ and
  $I_{A^{-1}B}$; note that $A^{-1}B$ and $AB^{-1}$ are conjugate.

  Since the Nielsen kernel $K(G)/G$ is homeomorphic to $\Delta/G$,
  it has only three boundary curves.  Looking at $F$, these must be the projections of
  the axes of $A,B$ and $A^{-1}B$ and all their conjugates.  Thus
  the boundary of the Nielsen region consists of the axes of the
  generators
  and the axes of all their conjugates and $F_K$ is a fundamental domain as
  claimed.
  \end{proof}
\subsection{Weierstrass points of the Schottky double $S$}
\label{sec:wps schottky double nonint}

 We now want to reflect the lines in $\partial\Delta$ so we
 introduce  subscripts to distinguish between the lines and their reflections
 as oriented line segments in $\oC$. Set
  $L_{int}=L$, $L_{ext}=R_{\partial \D}(L)$, etc.  The lines and
  their reflections in $\partial\D$ form
 triple of circles which we
denote by $(\C_L, \C_{L_A},\C_{L_B})$.  The circles are oriented
so that $\C_L = L_{int} \cup (L_{ext})^{-1}$.  Note further that
these three circles also bound a domain in $\oC$ and in fact, a
domain in $\Omega(G)$.

\begin{thm}
\label{thm:wps of Schottky double nonint} Suppose $G=\langle
C,D\rangle$ is a discrete free Fuchsian group such that the axes
of $C$ and $D$ do not intersect. Let $A,B$ be the stopping
generators determined by the Gilman-Maskit algorithm and let the
points $p_L,q_L$,$p_{L_B},q_{L_B}$,$p_{L_A},q_{L_A}$ be the
respective intersections of $L,L_A,L_B$ with $\partial\Delta$
(labelled in clockwise order beginning with $p_L$ and ending with
$q_{L_A}$) and so that
$(q_L,p_{L_B})$,$(q_{L_B},p_{L_A})$,$(q_{L_A},p_L)$ are disjoint
segments of $\partial\Delta$ respectively in the intervals of
discontinuity bounded by the axes $Ax_B$, $Ax_{A^{-1}B}$ and
$Ax_A$. Then the Schottky double $S$ is a compact Riemann surface
of genus two and the projections of these points onto $S$ under
the projection $\Omega(G) \rightarrow \Omega(G)/G=S$ are the
Weierstrass points of $S$.
\end{thm}

\begin{proof}  To simplify notation, we denote the
reflections $R_{\C_L}$ by $R_L$, $R_{\C_{L_A}}$ by $R_{L_B}$ and
$R_{\C_{L_B}}$ by $R_{L_B}$. Consider the circles
$\C_{L_{\bar{A}}}= R_{L}(\C_{L_A})=A(\C_{L_A})$ and
$\C_{L_{\bar{B}}}=\C_{R_{L}(L_B)}=B(\C_{L_B})$.

Since $(L,L_A,L_B)$ bound a domain in $\D$, the four circles
$\C_{L_A}$, $\C_{L_{\bar{A}}}$, $\C_{L_B}$,  and $
\C_{L_{\bar{B}}}$ bound a  domain $F$ in $\Omega$.

 This domain is the union of
the domain $F$ of figure~\ref{fig:s7} and its reflection in
$\partial\D$.  Moreover, $A$ maps the exterior of $\C_{L_A}$ to
the interior of $\C_{L_{\bar{A}}}$ and $B$ maps the exterior of
$\C_{L_B}$ to the interior of $\C_{L_{\bar{B}}}$. One can easily
verify that $F$ is a fundamental domain for $G=\langle A,B
\rangle$ and that $\Omega(G)/G$ is a compact Riemann surface of
genus two.

Let $R_{\partial\Delta}$ denote reflection in $\partial\Delta$ so
that $R_{\partial\Delta}: \Di \rightarrow \De$.  Set
$E=R_{\partial\Delta}R_L$.  The product of reflections in a pair
of intersecting circles is a rotation (an elliptic element of
$\MM$) about the intersection points of the circles (the fixed
points) with rotation angle equal to twice the angle between the
circles. As $E$ is the elliptic of order two  with fixed points
$p_L,q_L$, it takes any circle through $p_L$ and $q_L$ into itself
and interchanges the segments between the endpoints. In
particular, it sends $\C_L$ to itself and sends $L_{int}$ to
$(L_{ext})^{-1}$.

One easily verifies that  $EAE^{-1}=A^{-1}$ and $EBE^{-1}=B^{-1}$
so that $E$ is an orientation preserving  conformal involution
from $\Omega(G)$ to itself that
  conjugates $G$ to itself.  It therefore induces a
conformal involution $j: S \rightarrow S$.

Similarly, $E_A=R_{\partial\Delta}R_{L_A}=EA$ and
$E_B=R_{\partial\Delta}R_{L_B}=EB$ are elliptic elements of order
two with fixed points $\{p_{L_A},q_{L_A}\}$ and
$\{p_{L_B},q_{L_B}\}$ respectively. Therefore,
$E(p_{L_A})=A(p_{L_A})$, $E(q_{L_A})=A(q_{L_A})$ and
$E(p_{L_B})=B(p_{L_B})$,$E(q_{L_B})=B(q_{L_B})$.  Since $E_AE^{-1}
\in G$ and $E_BE^{-1} \in G$, we conclude that $E_A$, $E_B$ and
$E$ induce   the same conformal involution $j$. Moreover,
  the projections of
$p_L,q_L,p_{L_A},q_{L_A},p_{L_B},q_{L_B}$ are distinct points
each fixed by $j$.  Since $j$ is a conformal involution  fixing
six distinct points, it is the hyperelliptic involution on $S$.
\end{proof}

The reflection $R_{\partial\Delta}$  also normalizes
 $G$
and thus descends to a homeomorphism of $S$ which is the
anti-conformal involution $J$ (see
section~\ref{sec:limitsetsetc}). As we saw in
section~\ref{sec:doubles}, the clockwise arc of $\partial{\Delta}$
joining $q_{L_A}$ to $A(q_{L_A})=R_L(q_{L_A})$ projects to
  the central curve  of the funnel on $S$ corresponding to the projection of
interval of discontinuity bounded by the axis of $A$. Similarly,
the clockwise arc from $R_L(p_{L_B})$ to $p_{L_B}$ projects to the
central curve  of the funnel on $S$ corresponding to the
projection of the interval of discontinuity bounded by the axis of
$B$.   In addition, the pair of segments joining
(counterclockwise) $q_{L_B}$ to $p_{L_A}$ and $R_L(p_{L_A})$ to
$R_L(q_{L_B})$ together project to the central curve
 of the third funnel corresponding to the
projection of the interval of discontinuity bounded by the axis of $A^{-1}B$. These
three simple curves are fixed point-wise by $J$ and they partition
$S$ into two spheres with three holes (pairs of pants).

The product $Jj$
 is again an anti-conformal
self-map of  $S$.  Its invariant curves are the projections of the
circles $\C_L, \C_{L_A}, \C_{L_B}$.  The homotopy classes of these
projections also partition $S$ into two spheres with three holes.

We remind the reader that $J$ is an anti-conformal map on $S$ that
fixes its three fixed curves point-wise and $j$ is a conformal map
that sends each fixed curve into itself as a point-set, but only
fixes two points on each curve and maps the curve into its
inverse.

\section{The convex hull and the convex core }
\label{sec:convexd}

 In this section we
 consider the group $G$ acting on $\PP$, $\PPi$ and $\PPe$ which
 contain respectively the convex hull  $\CCC(G)$ in $\IH$ and its two
 boundaries $\CCC_{int}$ and $\CCC_{ext}$. All the axes of
 elements of $G$ lie in $\PP$.

Note that the reflection $R_{\PP}$ interchanges $\PPi$ and $\PPe$.
Since $\PP \subset \IH$, in addition to factoring elements of $G$
as products of reflections in hyperbolic planes we can also factor
them  products of half turns about hyperbolic lines.  (See
\cite{Fench}).  If $M$ is a hyperbolic line in $\IH$, we denote
the half turn about $M$ by $H_M$.

\subsection{Fundamental domains in $\PP$ and $\CCC$}
\label{sec:funddomains in h3}

 Since  the axes of the generators $C$ and $D$ do not intersect and
 these axes lie  in the plane $\PP$
 (which is another model for $\D$), we can
apply the Gilman-Maskit algorithm to obtain stopping generators
$A,B$ as above. We again have a triple of lines $L,L_A,L_B$
defined as follows: $L$ is the hyperbolic line mutually orthogonal
to the axes of $A$ and $B$; it now lies in $\PP$. The hyperbolic
line $L_A$ in $\PP$ is determined by the factorization
$A=H_LH_{L_A}$ Similarly $L_B$ is the line in $\PP$ such that
$B=H_{L}H_{L_B}$. The condition for $A,B$ to be stopping
generators is that the lines $L, L_A, L_B$ bound a domain in
$\PP$. Set $L_{\bar{A}}=H_L(L_{A})$, $L_{\bar{B}}=H_L(L_{B})$.  It
is easy to check that $L_A$ is the mutual perpendicular to the
axes of $A$ and $A^{-1}B$ while $L_B$ is the mutual perpendicular
to the axes of $B$ and $A^{-1}B$

\begin{figure}[hbt]\centering
\psfrag{LbB}[cc]{$L_{\bar{B}}$}
 \psfrag{LbA}[cc]{$L_{\bar{A}}$}
\psfrag{PL}[cc]{$p_L$}
 \psfrag{L}[cc]{$L$}
\psfrag{QLA}[cc]{$q_{L_A}$}
 \psfrag{LA}[cc]{$L_A$}
\psfrag{PLA}[cc]{$p_{L_A}$}
 \psfrag{QLB}[cc]{$q_{L_B}$}
\psfrag{LB}[cc]{$L_{B}$}
 \psfrag{QL}[cc]{$L_B$}
\psfrag{PLB}[cc]{$p_{L_B}$}
 \psfrag{Delta}[cc]{$\Delta$}
\includegraphics[width=5in]{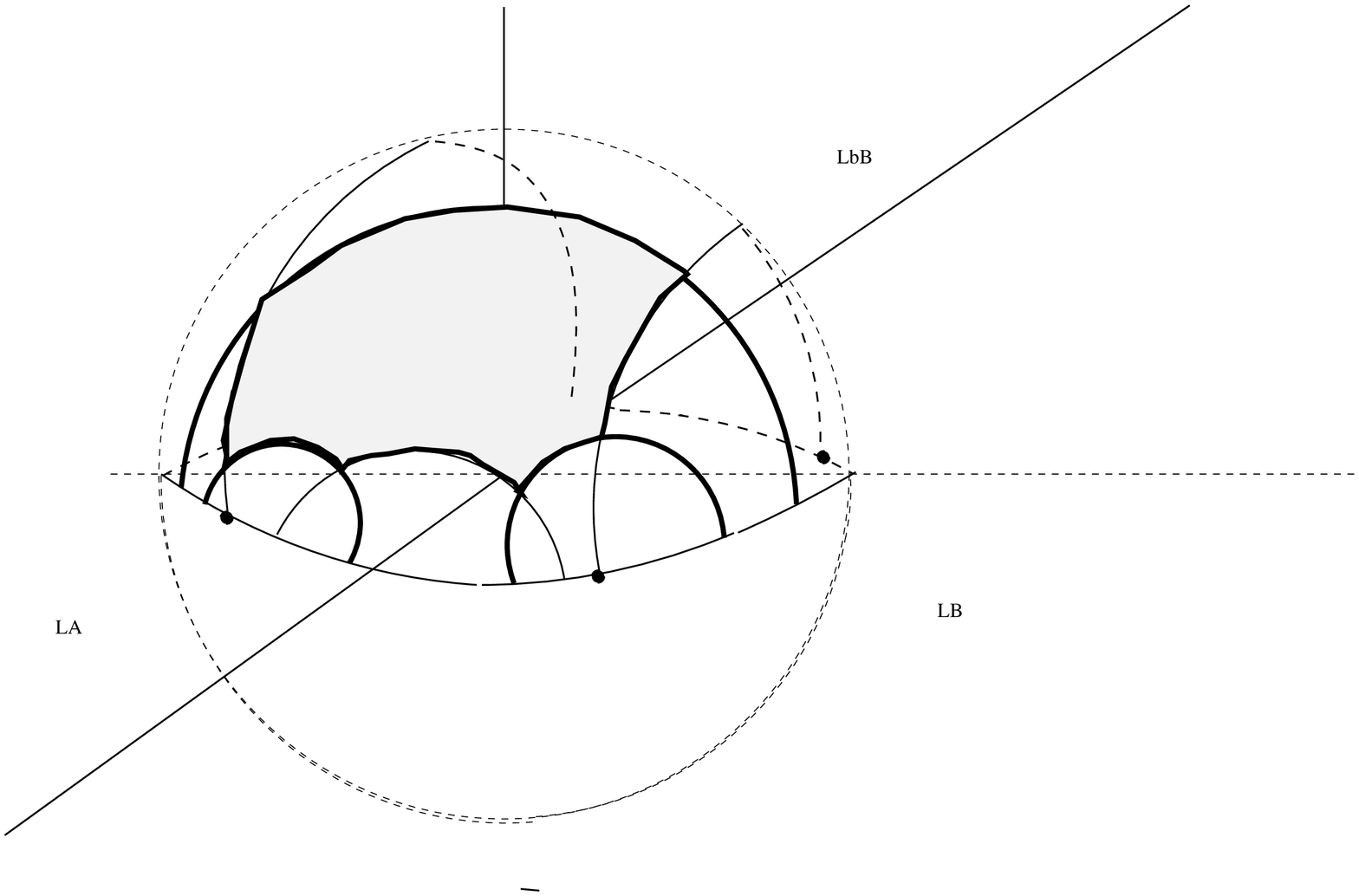}
\caption{The domain $\HHH \subset \PP \subset \IH$.
\label{figure:newf}}
\end{figure}

 Let $\HHH$ be the closed right-angled hexagon formed by
the six hyperbolic lines, $L,L_A,L_B$ and $Ax_A,Ax_B,Ax_{A^{-1}B}$
and shown in figure~\ref{figure:newf}.

\begin{prop}\label{prop:funddomaininP nonint} The interior of the
domain $F_{\CCC} =int(\HHH \cup H_{L}(\HHH))$ is a
fundamental domain for $G$ acting on the convex hull $\CCC(G)$.
\end{prop}

\begin{proof} Since $\CCC(G) \subset \PP$, it is same as the Nielsen region
$K(G_{\PP})$, and the proposition follows directly from
proposition~\ref{prop:funddomain nonint}.
\end{proof}

We define the domains  $F_{int} = F_{\CCC} \cap \PPi $ and
$F_{ext} = F_{\CCC} \cap \PPe$ respectively. Note that $F_{int}
\subset \CCC_{int} \subset \PPi$ and $F_{ext} \subset \CCC_{ext}
\subset \PPe$. An immediate corollary of the above proposition is

\begin{cor}
\label{cor:Fis corebdy} The domain $F_{int} \cup F_{ext}$ is a
fundamental domain for $G$ acting on $\CCC_{int} \cup
\CCC_{ext}=\partial\CCC$.
\end{cor}

If we identify the appropriate sides of $F_{\CCC}$ under the action of
 $A$ and $B$,
we obtain a sphere with three boundary curves that are geodesic in
the intrinsic metric. Therefore, identifying the sides of $F_{int}
\cup F_{ext}$ we obtain two such spheres.  Gluing their boundary
curves yields a surface of genus two as the boundary of the convex
core that  is isometric with the Nielsen double of $G$.

\section{The Pleated Surface  }
\label{sec:pleated surface -nonint} In this section we explicitly
construct the pleated surface $(\tilde{S},pl)$ where the
hyperbolic surface $\tilde{S}=S_K$. As usual we assume we have
simple stopping generators with disjoint axes.

To define the pleated surface we need first to construct a
hyperbolic surface of genus two.  We do this by constructing a
Fuchsian group and forming the quotient. The hyperbolic surface we
obtain is the Nielsen double, $S_K$.

\subsection{The Fuchsian group for the Nielsen double }
\label{sec:group gamma non int}
 We begin with the fundamental domain $F_K=\HHH \cup
 R_L(\HHH )\subset \D$
defined in section~\ref{sec:schottkynonint} for the Nielsen
kernel.
  Set   $F_{\Gamma} = F_K \cup
R_{Ax_A}(F_K) \subset \Delta$. Note that $F_{\Gamma}$ contains
four copies of the hexagon $\HHH$.
  Define a new group acting on $\Delta$ by $\Gamma=\langle
  A,A',B,B'\rangle$,
where $A'=R_{Ax_{A^{-1}B}}R_{Ax_{A}}$ and
$B'=R_{Ax_{B}}R_{Ax_{A}}$.

\begin{figure}[hbt]\centering
\psfrag{Rone}[cc]{$R_{Ax_A}(L(\HHH))$}
\psfrag{Rtwo}[cc]{$R_{Ax_A}(\HHH)$}
 \psfrag{PA}[cc]{$p_A$}
\psfrag{PB}[cc]{$p_B$}
 \psfrag{P}[cc]{$p$}
\psfrag{AxA}[cc]{$Ax_A$}
 \psfrag{LbB}[cc]{$L_{\bar{B}}$}
\psfrag{L\\bB}[cc]{$L_{\bar{B}}$}
\psfrag{LbA}[cc]{$L_{\bar{A}}$}
 \psfrag{PL}[cc]{$p_L$}
\psfrag{L}[cc]{$L$}
\psfrag{QLA}[cc]{$q_{L_A}$}
\psfrag{LA}[cc]{$L_A$}
 \psfrag{PLA}[cc]{$p_{L_A}$}
\psfrag{QLB}[cc]{$q_{L_B}$}
\psfrag{LB}[cc]{$L_{B}$}
\psfrag{QL}[cc]{$L_B$}
 \psfrag{PLB}[cc]{$p_{L_B}$}
\psfrag{\Delta}[cc]{$\;$}

\psfrag{LHS}[cc]{$L(\HHH)$}
 \psfrag{LH}[cc]{$L(\HHH)$}
\psfrag{R_{Ax_A}(L(H))}[cc]{$R_{Ax_A}(L(\HHH))$}
\psfrag{R_{Ax_A}(H)}[cc]{$R_{Ax_A}(\HHH)$}
 \psfrag{MA}[cc]{$M_A$}
\psfrag{MB}[cc]{$M_B$}
 \psfrag{MBA}[cc]{$M_{\bar{A}}$}
 \psfrag{MbB}[cc]{$M_{\bar{B}}$}
 \psfrag{M'bA}[cc]{$M'_{\bar{A}}$}
 \psfrag{MB}[cc]{$M_{B}$}
 \psfrag{MbA}[cc]{$M_{\bar{A}}$}
 \psfrag{MA}[cc]{$M_{A}$}
 \psfrag{H}[cc]{$\HHH$}
   \psfrag{AXA}[cc]{$Ax_A$}
   \psfrag{AXB}[cc]{$Ax_B$}
\includegraphics[width=5in]{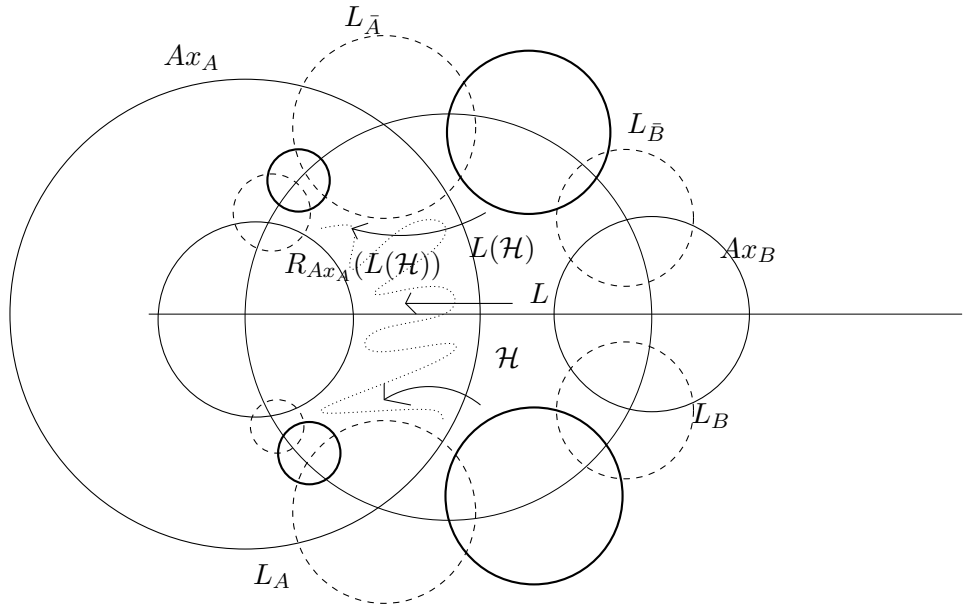}
\caption{The fundamental domain  for $\Gamma$ is a union of four
copies of $\HHH$. The circle sides, $L_A,L_{\bar{A}}$,
$L_B,L_{\bar{B}}$ and their reflections in $Ax_A$ are drawn
dotted. The axis sides are drawn solid. \label{figure:srefl}}
\end{figure}

\begin{prop}
\label{prop:genus2funddomain} The domain  $F_{\Gamma}$ is a
fundamental domain for the group $\Gamma=\langle A,A',B,B'\rangle$
acting on $\Delta$. The Fuchsian surface $\tilde{S}=\Delta/\Gamma$
is a Riemann surface of genus two isometric to the Nielsen double
$S_K$.
\end{prop}

\begin{proof}
It suffices to show that the sides are identified in pairs by
elements of $\Gamma$ and the angles at the vertices satisfy the
conditions of the Poincar\'e polygon theorem (\cite{Mas}).

In $F_K$ the pair of sides that are segments of $L_B$ and
$L_{\bar{B}}$ are identified by $B$; similarly the sides of $L_A$
and $L_{\bar{A}}$ in $F_{\Gamma}$ are identified by $A$. In
$F_{\Gamma}$ the pair of sides that are segments of $Ax_B$ and
$R_{Ax_A}(Ax_B)$ are identified by $(B')^{-1}$ since $R_{Ax_A}:R_{Ax_A}(Ax_B) \rightarrow Ax_B$ and $R_{Ax_B}:Ax_B \rightarrow Ax_B$.  Then the sides that are
segments of $R_{Ax_{A}}(L_B)$ and $R_{Ax_{A}}(L_{\bar{B}})$ are
identified by $B'B(B')^{-1}$.  Similarly, The sides that are segments of
$Ax_{A^{-1}B}$ and $R_{Ax_{A}}(Ax_{A^{-1}B})$ are identified by
$(A')^{-1}$ and the sides that are segments of  $Ax_{AB^{-1}}$ and
$R_{Ax_{A}}(Ax_{AB^{-1}})$ are identified by $AA'A^{-1}$. (See
figure~\ref{figure:srefl}.)

It is easy to check that the twelve  vertices of $F_{\Gamma}$ are
partitioned into three cycles of four vertices. Denote the three
vertices of $\HHH$ in $F_{\Gamma}$ as follows: let $q_A$ be the
intersection of $Ax_{A^{-1}B}$ with $L_A$, let $p_B$ be the
intersection of $Ax_{A^{-1}B}$ with $L_B$ and let $q_B$ be the
intersection of $Ax_B$ and $L_B$.  The cycles are then
$$\{p_A, A'(p_A), AA'(p_A), A(p_A)\},$$
$$\{p_B,B(p_B),B'B(p_B),B'(p_B)\},$$
$$\{q_B,A'(q_B),B(q_B),B'B(B')^{-1}A'(q_B)\}$$ As each angle in
$F_{\Gamma}$ is a right angle, the sum of the  angles at each
cycle of vertices is $2\pi$.

We conclude that $F_{\Gamma}$ is a fundamental domain for the
group generated by the side pairing transformations. Since these
are generated by the generators of $\Gamma$, $F_{\Gamma}$ is a
fundamental domain for $\Gamma$.

Sides of  $F_{\Gamma}$ either are termed axis sides or reflection
sides depending upon whether they are segments of an axis of an
element of $G$  or segments of a reflection line of $G$.

Now $F_K$ is composed of two copies of $\HHH$, and, as we saw in
section~\ref{sec:schottkynonint},
 identifying  non-axis sides (segments of $L_A,L_{\bar A},L_B,L_{\bar B}$)
 yields a sphere with three holes. The boundaries of the holes are
 the identified axes $Ax_A,Ax_B,Ax_{AB^{-1}}$.
Similarly $R_{Ax_A}(F_K)$ with non-axis sides identified is
another three holed sphere. The axis sides
 $(Ax_A,Ax_B,Ax_{AB^{-1}})$ of $F_{\Gamma}$ are identified as indicated above,
 and the two three holed
spheres join up to form a compact surface $\tilde{S}$ of genus
two.  Since $F_{\Gamma}$ is just $F_K$ doubled, the surface
$\tilde{S}$ is isometric to the Nielsen double $S_K$.
\end{proof}

\subsection{The pleating map }
\label{sec:pleating map nonint}
 To construct the  pleating map
$pl: \D/\Gamma \rightarrow \IH/G$, we begin by defining a pleating
map $PL:\D \rightarrow \IH$. We first define $PL$  on $F_{\Gamma}$
and then extend by the groups $\Gamma$ and $G$.
 Recall that the
domain $F_{\Gamma} \subset \D$ is a union of two copies of $F_K$,
and $F_{int} \subset \PP_{int}$ and $F_{ext} \subset \PP_{ext}$ as
defined in section~\ref{sec:funddomains in h3}  are each isometric
to $F_K$.
 Set
$PL: F_K \rightarrow F_{ext}$. Next  set $PL: R_{Ax_A}(F_K)
\rightarrow F_{int}$.
 Define the group homomorphism $\phi: \Gamma \rightarrow G$ by first
defining  it on the generators:
$$A \mapsto A, B \mapsto B, A' \mapsto  id, B' \mapsto id.$$
 To show that this map on generators gives a group homomorphism
from $\Gamma$ to $G$ we must show that $\phi$ preserves the
defining relation(s)
 of $\Gamma$. We know from the Poincar\'e polygon theorem
  that the reflection relations and the cycle relations
form  a complete set of relations for $\Gamma$  There are no
reflection relations and it is easy to calculate that substituting
$\phi(A)$, $\phi(A')$, $\phi(B)$,
 and $ \phi(B')$ into the cycle relation gives the identity.

We now extend $PL$ to a map from $\D$ into $\IH$. That is, for $x
\in \D$ let $g \in \Gamma$ be chosen so that $g(x) \in F_{\Gamma}$
and set $PL(x) = \phi(g^{-1}) \circ PL \circ g(x)$. This is well
defined  when  $g(x)$ lies interior to $F_{\Gamma}$ and when
$g(x)$ lies on a boundary curve of $F_{\Gamma}$ or on $Ax_A$ it is
defined by continuity. Thus,  $(\D,PL)$ is a pleated surface with
image in $\IH$ whose pleating locus is the image of the axis sides
of $F_{\Gamma}$ and all their images under $\Gamma$.

We now define the map $pl:\tilde{S} \rightarrow \IH/G$ by taking
quotients.   The map $pl$ is clearly a hyperbolic isometry since
$PL$ is and its pleating locus consists of the images of  $Ax_A$,
$Ax_B$ and $Ax_{A^{-1}B}$.

By corollary~\ref{cor:Fis corebdy}, we can apply the group $G$ to
identify $PL(F_{\Gamma})$ with $\partial\CCC$, and taking
quotients, identify $pl(S_K)$ with the convex core boundary
$\partial{N}=\CCC_{int}/G\cup \CCC_{ext}/G$.

\subsection{The Weierstrass points  of the Nielsen double $S_K$ }
\label{sec:wps nielsen double nonint}

 We now determine the
Weierstrass points of $S_K$. This is very reminiscent of our
earlier discussion of the Weierstrass points of the Schottky
double, except that now, instead of using the simple reflection in
$\partial\Delta$ we need to use reflections in the axis sides of
$\HHH$ that correspond to the generating triple with simple axes.

\begin{quote}
  Let $E_0=R_LR_{Ax_A}$. Let $p$ be the fixed point of $E_0$ .\\
Let $E_1=R_LR_{Ax_B}$. Let $q$ be the fixed point of $E_1$\\
Let $E_2=R_{L_A}R_{Ax_A}$. Let $q_A$ be the fixed point of $E_2$.\\
\end{quote}

Note that $E_1=B'E_0$ and $E_2=A^{-1}E_0$.  The vertices of
$F_{\Gamma}$ are also fixed points of order two  elliptics: $p_A$
is the fixed point of $E_3=E_2A'^{-1}=A^{-1}E_0A'^{-1}$, $q_B$ is
the fixed point of $E_4=A^{-1}E_0A'B^{-1}A=E_3B^{-1}A$ and $p_B$
is the fixed point of $E_5=E_1B=B'E_0B$.

\begin{thm}
\label{thm:wps of Nielsen double nonint}The Weierstrass points of
$S_K$ are the projections of the six points $p,q,p_A,q_A,p_B,q_B$.
\end{thm}
\begin{proof}
Since $E_0=R_LR_{Ax_A}$, we see that
$E_0(F_{\Gamma})=F_{\Gamma}$ and verify
 that
$E_0$ normalizes $\Gamma$.
 It follows that $E_0$  induces a conformal involution $\tilde{j}$ on
$\tilde{S}$.

Similarly, $E_i$, $i=1,\ldots 5$ preserves the tiling by images of
$F_{\Gamma}$ and thus  normalizes $\Gamma$ so that $E_i$ projects to
a conformal involution. Since the products of pairs of the $E_i$'s
are elements of $\Gamma$, they all project to the same conformal
involution $\tilde{j}$ on $\tilde{S}$. As the six points are
inequivalent under $\Gamma$, they project to six distinct fixed
points of the involution $\tilde{j}$, characterizing it as the
hyperelliptic involution.
\end{proof}

\section{Weierstrass points and lines in the handlebody}
\label{sec:gen weierstrass non int} In this section we give an
explicit description of the handlebody $H$ as $\cal{S} \times I$
where $\cal{S}$ is a compact surface of genus two and $I$ is the
interval $[0,\infty]$. We show that $H$ admits a unique order two
isometry $j$ that we call the {\em hyperelliptic isometry of $H$}.
It fixes six unique geodesic line segments, which we call the {\sl
Weierstrass lines} of $H$. The Weierstrass lines are, of course,
the analogue of the Weierstrass points on the boundary surface.

 In addition, we show that $H$ is
foliated by surfaces $\cal{S}(s)=(\cal{S},s)$ that are at
hyperbolic distance $s$ from the convex core, and that each
$\cal{S}(s)$ is fixed by $j$. The restriction of $j$ to
$\cal{S}(s)$ fixes six {\sl generalized Weierstrass points} on
$\cal{S}(s)$. On each $\cal{S}(s)$  we also
 find an orientation reversing involution $J$ that
fixes curves through the  generalized Weierstrass points.

\subsection{Construction of the foliation}
\label{sec:equisurface disjoint}

We construct the {\sl hyperelliptic isometry for $H$} as follows.
We saw, in section~\ref{sec:funddomains in h3}, that $G$ is
generated by products of even numbers of half-turns in the lines
$(L,L_A,L_B)$. It follows that   $G$ is a normal subgroup of index
two in the group $ \langle H_{L}, H_{L_A}, H_{L_B} \rangle$. Since
$ H_{L}$, $H_{L_A}$ and $H_{L_B}$ all differ by elements of $G$,
their actions by conjugation on $G$ all induce the same order two
automorphism $j$ under the projection $\pi:\IH \rightarrow
\IH/G=H$. The fixed points of $j$ in $H$ are precisely the points
on the images of the lines $ L, L_A,{L_B}$; denote the projection
of the lines by $(\pi(L), \pi(L_A),\pi(L_B))$.

We work with the lines $L$ and $\pi(L)$ but we have analogous
statements for the  lines $L_A$ and $L_B$ and their projections.
As we saw in theorem~\ref{thm:wps of Schottky double nonint}, the
endpoints $(p,q)$ of $L$ on $\partial\D$ project to  Weierstrass
points on the Schottky double $S=\partial{H}=\cal{S}(\infty)$.

 In theorem~\ref{thm:wps of Nielsen double nonint} we
found the Weierstrass points of the Nielsen double $S_K$. In the
fundamental domain $F_{\Gamma} \subset \D$ for $\Gamma$ we found
points $p \in F_{\Gamma}$ and $q \in
\partial{F_{\Gamma}}$ such that $PL(p)$ and $PL(q)$ both lie on the
line $L \in \PP$.  For the moment, call the  hyperelliptic
involution on $S_K$, $j_{S_K}$ and the hyperelliptic involution on
$H$, $j_H$. By construction, the pleating map $pl$ commutes with
these involutions: $pl \circ j_{S_K}=j_H \circ pl$.  It follows
that the projections $\tilde{p}=\pi(PL(p))$ and
$\tilde{q}=\pi(PL(q))$ are fixed under $j_H$. By the same
argument, we find two other  pairs of points
$\pi(PL(p_A)),\pi(PL(q_A))$ on $\pi(L_A)$ and
$\pi(PL(p_B)),\pi(PL(q_B))$ on $\pi(L_B)$ fixed by $j_H$.  They
lie on the genus two surface $\partial{N}=\cal{S}(0)$ and are its
{\em generalized Weierstrass points}.

The segment of the line $\pi(L)$ from $\tilde{p}$ to $\tilde{q}$
lies inside $N$ and its complement consists of two segments: one
from $\tilde{p}$ to $\pi(p) \in \cal{S}(\infty)$ which we denote
$L_p$, and the other from $\tilde{q}$ to $\pi(q) \in
\cal{S}(\infty)$ which we denote by $L_q$. Each of these is
point-wise fixed under $j_H$. We call the two segments outside
$N$, {\em the generalized Weierstrass lines of H}. Similarly, we
have two other pairs of generalized Weierstrass lines
$(\pi(L_{A_p}),\pi(L_{A_q}))$ and $(\pi(L_{B_p}),\pi(L_{B_q}))$.

 We can parameterize these lines by hyperbolic length. We next
 construct the family of surfaces $\cal{S}(s)$. We will show that
 $j_H$, which we denote again simply by $j$ since it will not
 cause confusion, is an order two isometry of $\cal{S}(s)$ that has as its set of fixed points
  the  intersection points of these six lines with the surface.

For a point $p \in \PP$, let $V(t)$, $t \in (-\infty,\infty)$ be
the line in $\IH$ perpendicular to $\PP$ at $p$ and parameterized
so that $s=|t|$ is hyperbolic arc length and oriented  such that
$V(s)$ is exterior to $\PP$ and $V(-s)$  is interior to $\PP$.
There is a family of surfaces $\Pi(\pm s)$ with boundary
$\partial\Delta$, passing through the point $V( \pm s)$ such that
the distance from any point on $\Pi( \pm s)$ to $\PP$ is $s$. We
 call them {\em equidistant surfaces}. We can think of $
\Pi(s) \cup \Pi(-s)$ as the (full) equidistant surface from $\PP$
at a distance $s$, with the points lying above $\PP$ having
directed distance $s$ and the points below $\PP$ having directed
distance $-s$. (See \cite{Fench} chapters III.4 and IV.5).  Note
that since $G$ acts on $\IH$ by isometries, it leaves the
hyper-surfaces $\Pi(\pm s)$ invariant.

\begin{figure}[hbt]\centering
\includegraphics[width=2in, height=2in]{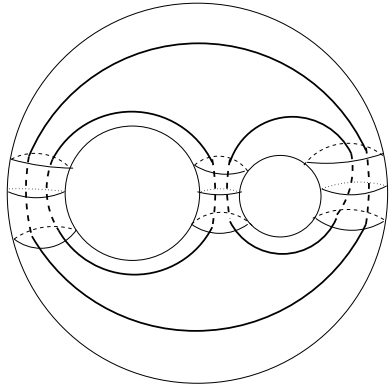}
\caption{The convex core inside the handlebody
\label{figure:surf1}}
\end{figure}


Using orthogonal projection, we can project ${\CCC_{ext}}$ onto
$\Pi(s)$ to obtain $\CCC_{ext}(s)$ and project $\CCC_{int}(s)$
onto $\Pi(-s)$  to obtain $\CCC_{int}(s)$.

Each quotient, $N_{ext}(s)= \CCC_{ext} /G$ and $N_{int}(s) =
\CCC_{int}(s)/G$, is topologically a sphere with three boundary
curves. In analogy with our construction of the Schottky double,
we want to join these boundary curves, in pairs, by  funnels to
form the surface $\cal{S}(s)$  (see figure~\ref{figure:surf1}).

Let  $Eq_A(s)$ be the equidistant cylinder in $\IH$ about $Ax_A$.
Note that, for each $s>0$,  $A$ maps $Eq_A(s)$ to itself. Because
$Ax_A$ lies in $\PP$, $A(\pm s)= Eq_A(s)\cap \Pi(\pm s)$ is the
orthogonal projection of $Ax_A$ onto $\Pi(\pm s)$; notice that
$Eq_A(s)$ and $\Pi(\pm s)$ are tangent along $A(\pm s)$.  The
curves $A(\pm s)$ intersect $\CCC_{ext}$ and $\CCC_{int}$ in
boundary curves. The curves $A(\pm s)$ divide $Eq_A(s)$ into two
pieces, each a doubly infinite topological strip. Both of these
strips intersect $\PP$, one on each side of the axis of $A$. Thus
one strip intersects $\CCC$ and one does not;  we denote the strip
disjoint from $\CCC$ by $Rec_A(s)$.

Similarly, we form $Rec_B(s)$, $Rec_{A^{-1}B}(s)$,
$Rec_{AB^{-1}}(s)$, \ldots, for each of the boundary curves of
$\CCC$. The quotient of each of these strips under the action of
$G$ is a  funnel; there are three distinct funnels and we denote
them $F_A(s), F_B(s)$ and $F_{A^{-1}B}(s)$.

The  surface $\cal{S}(s)$ is defined as
$$\cal{S}(s)=N_{ext}(s) \cup N_{int}(s) \cup F_A(s)\cup F_B(s) \cup F_{A^{-1}B}(s) $$
 We have
\begin{prop} \label{prop:calS is invariant disj} The surface
$\cal{S}(s)$ is invariant under the involution $j$. \end{prop}
\begin{proof} By construction, the distance from each point of $\cal{S}(s)$ to
the convex core $N$ is $s$.  As the involution $j:H \rightarrow H$
is an isometry and it preserves $N$, it leaves $\cal{S}(s)$
invariant. \end{proof}

We want to find the fixed points of $j$ acting on $\cal{S}(s)$. To
this end consider the intersection points of the line segments
$L_p(s)$ and $L_{A_p}(s)$ with $Rec_A(s)$;  denote them
respectively by $p(s)$ and $p_A(s)$.  Similarly, define the points
$q(s)$ and $q_B(s)$ as the intersection points of the line
segments $L_q(s)$ and $L_{B_q}(s)$ with $Rec_B(s)$ and the define
the points $q_A(s)$ and $p_B(s)$ as the intersection points of the
line segments $L_{A_q}(s)$ and $L_{B_p}(s)$ with
$Rec_{A^{-1}B}(s)$.

\begin{prop} The projections of the points $p(s),q(s)$,
$p_A(s),q_A(s)$ and  $p_B(s),q_B(s)$ lie on the surface
$\cal{S}(s) \subset H$ and comprise the set of fixed points of $j$
restricted to $\cal{S}(s)$.
\end{prop}
\begin{proof} Each of these points lies on one of the six
generalized Weierstrass lines in $H$ at distance $s$ from the
convex core. These lines are distinct so the points are distinct.
Since these lines are fixed pointwise by $j$, and since $j$ leaves
$\cal{S}(s)$ invariant, these are fixed points of $j$ on
$\cal{S}(s)$ as claimed.
\end{proof}

We call these projected points the {\em generalized Weierstrass
points} of $\cal{S}(s)$.

We summarize these results as
\begin{thm}
\label{thm:foliation nonint} The projections of the
 lines $L,L_A,L_B$ in $ \PP \subset H$ are point-wise fixed by the involution
 $j$.  The projection of each line  consists
of three segments. The endpoints of the internal segment, lying in
the convex core,  are  generalized Weierstrass points of
$(\cal{S},0)$. Each line has two  external segments, parameterized
by hyperbolic arc length and called generalized Weierstrass lines.
The points at distance $s$ on the generalized Weierstrass lines
are the generalized Weierstrass points of the surface
$(\cal{S},s)$. The endpoints of the external segments lie on the
boundary the handlebody and are the Weierstrass points of the
boundary surface $S=(\cal{S},\infty)$.
\end{thm}

\subsection{Anticonformal involutions}
\label{sec:anticonf nonint} We can also construct the
anti-conformal involution $J$ acting on $H$ by defining it to be
the self-map of $H$ induced by the anti-conformal reflection in
the plane $\PP$, $R_{\PP}$. Since the axes of the generators $A$
and $B$ of $G$ lie in $\PP$, $R_{\PP}$ fixes theses axes
point-wise and since $R_{\PP}$ is an orientation reversing map, it
induces an orientation reversing map on $H$.

We want to see how $J$ acts on
$$\cal{S}(s)=N_{ext}(s) \cup N_{int}(s) \cup F_A(s)\cup F_B(s) \cup F_{A^{-1}B}(s) $$
Since $R_{\PP}$ is an isometry,
$R_{\PP}(\CCC_{ext}(s))=\CCC_{int}(s)$ and $J(N_{ext})=N_{int}$.
 The involution maps $Rec_A(s)$ to itself, fixing the curve that
is the intersection of the plane $\PP$ with $Rec_A(s)$. We call
the projection of this curve the central curve of the funnel.
Since $R_{\PP}$ interchanges the part of $Rec_A(s)$ on one side of
$\PP$ with the part on the other, we deduce that $J$ maps the
 funnel $F_A(s)$ to itself, interchanging its boundary
curves and fixing the central curve.

There is another  anti-conformal self-map  of $H$ induced by
$R_{\PP} \circ H_L$; $\hat{J}=J \circ j$.  To see how $\hat{J}$
acts on $\cal{S}(s)$ we first look at the fundamental domains
$F_{ext}(s)$ and $F_{int}(s)$.  Both $j$ and $J$ interchange these
so $\hat{J}$ leaves each of them invariant.  We look at
$F_{ext}(s)$: Let $\PP_L$ be the plane orthogonal to $\PP$ through
the line $L$ and let $\PP_L(s)$ be the intersection of $\PP_L$ and
$\Pi(s)$; equivalently, $\PP_L(s)$ is the orthogonal projection of
$L$ onto $\Pi(s)$.  Looking back at the construction of $F_{ext}$,
we see that it is symmetric about $L$ and the symmetry is given by
$R_{\PP} \circ H_L$; thus $F_{ext}(s)$ is symmetric about
$\PP_L(s)$ by the same map.  It follows that $\hat{J}$ maps
$N_{ext}(s)$ to itself with fixed curve the projection of
$\PP_L(s)$.  We have an analogous symmetry for $N_{int}(s)$.

The plane $\PP_L$ also bisects the region $Rec_A(s)$. The map
$R_{\PP} \circ H_L$ interchanges the two halves of $Rec_A(s)$.
Thus $\hat{J}$ maps each  funnel to itself; it leaves each
boundary curve and the central curve invariant, but changes its
orientation. On the central curve, it leaves the Weierstrass
points fixed but interchanges the segments they divide the curve
into.
\nopagebreak
\part{INTERSECTING AXES: $\T [C,D] < - 2$}
\nopagebreak
\label{part:int}  In this part we always assume the axes of the
generators of $G$ intersect.  The organization  parallels part 1.
We begin by finding fundamental domains for $G$ acting on $\D$ and
on $K(G)$. Next, we find the Weierstrass points for $S$.  We then
turn to $\IH$ and find fundamental domains for $G$ acting on $\PP$
and on $\CCC$. We construct the Fuchsian group uniformizing $S_K$
and the pleating map sending $K(G)$ onto $\partial{\CCC}$. We then
use the pleating map to find the Weierstrass points of $S_K$.
Finally, we work with the handlebody $H$ and construct the
Weierstrass lines and the equidistant surfaces.  We conclude with
a discussion of anti-conformal involutions on the various
quotients.

\section{Subspaces of $\Delta$ and quotients by $G$}
\label{sec:schottky int}  Any pair of generators in the
intersecting axes case are geometric (see \cite{KACT},
\cite{KANN}).
 To unify our
discussion with the non-intersecting case and to simplify our
constructions and proofs below, however, we again replace the
generators by some Nielsen equivalent generators with special
properties.

Let $E_x$ denote the elliptic element of order two that leaves the
disc $\Delta$ invariant and fixes the point $x$. We let $p= Ax_C
\cap Ax_D$ and choose $p_C$ and $p_D$  so that $C= E_p \circ
E_{p_C}$ and $D = E_p \circ E_{p_D}$.  Note that $p$ and $p_C$
both lie on $Ax_C$ and that $\rho_{\D}(p,p_C)$ is half the
translation length of $C$, where $\rho_{\D}$ is the usual
hyperbolic metric in $\D$.

In a manner analogous to  the steps of the algorithm in the case
when the commutator is elliptic (\cite{Galg}), one can replace the
original generators by  generators $A,B$ so that the distance
between the three points $p,p_A,p_B$ is minimal (among all Nielsen
equivalent pairs of generators).   That is, $\rho_{\D}(p,p_A) \le
\rho_{\D}(p, p_B) \le \rho_{\D}(p_A,p_B)$ and the traces of $A$
and $B$ are minimal among all pairs of generators. Geometrically,
this condition says that the triangle with vertices $p,p_A,p_B$ is
acute.  This acute triangle is the analog of the domain bounded by
the lines $L,L_A,L_B$.

 By analogy with
the disjoint axes case, we will call the minimal trace generators
{\em stopping generators}.
   In this part we
 will always assume that the generators $A$ and $B$ are the
 stopping generators of $G$.

 \subsection{Fundamental Domains in $\D$
and $K(G)$} \label{sec:funddomains int}

We begin by constructing a fundamental domain for $G$ acting on
$\D$ whose side pairings are the stopping generators.

 We observe that $Ax_{[A,B]}$ is disjoint
from the axes $Ax_A$ and $Ax_B$ as follows. Following \cite{Mat}
construct $h$, the perpendicular from $p_B$ to
 $Ax_A$. By the trace minimality of the stopping generators,
 $h$ intersects $Ax_A$ between
 $p_A$ and $A^{-1}(p)$ (or vice-versa) \cite{Galg}.
Let $\delta$ be a geodesic perpendicular to $h$ passing through
$p_B$. The transformation $E_p\circ E_{p_A}\circ E_{p_B}$
identifies $\delta$ and $A(\delta)$. It is easily seen to be the
square root of the commutator $[A,B^{-1}]$ so that both have the
same axis, the common perpendicular to $\delta$ and $A(\delta)$.
 By
construction,  this axis must be  disjoint from $Ax_A$ and $Ax_B$.
The construction also assures that the axes of the commutators
$[B^{-1},A^{-1}]$,$[A,B^{-1}]$,$[B,A]$,$[A^{-1},B]$ are disjoint
from the axes of $A$ and $B$.

\begin{figure}[hbt]\centering
\psfrag{PA&}[cc]{$p_A$}
 \psfrag{PB&}[cc]{$p_B$}
\psfrag{P&}[cc]{$p$}

 \psfrag{AxA}[cc]{$Ax_A$}

 \psfrag{AxBA}[cc]{$Ax_{[A,B]}$}
\psfrag{AxA1BA}[cc]{$Ax_{[A^{-1},B]}$}
\psfrag{AxB1A1}[cc]{$Ax_{[B^{-1},A^{-1}]}$}
\psfrag{AxAB1}[cc]{$Ax_{[A,B^{-1}]}$}
\psfrag{LbB}[cc]{$L_{\bar{B}}$}
 \psfrag{LbA}[cc]{$L_{\bar{A}}$}
\psfrag{PL}[cc]{$p_L$}
 \psfrag{L}[cc]{$L$}
\psfrag{QLA}[cc]{$q_{L_A}$}
 \psfrag{LA}[cc]{$L_A$}
\psfrag{PLA}[cc]{$p_{L_A}$}
 \psfrag{QLB}[cc]{$q_{L_B}$}
\psfrag{LB}[cc]{$L_{\bar{B}}$}
 \psfrag{QL}[cc]{$L_B$}
\psfrag{PLB}[cc]{$p_{L_B}$}
 \psfrag{\Delta}[cc]{$\Delta$}
\psfrag{Ax[A,B^{-1}]}[cc]{$Ax_{[A,B^{-1}]}$}
\psfrag{Ax[B,A]}[cc]{$Ax_{[B,A]}$}
\psfrag{Ax[B^{-1},A^{-1}]}[cc]{$Ax_{[B^{-1},A^{-1}]}$}
\psfrag{Ax[A^{-1},BA]}[cc]{$Ax_{[A^{-1},BA]}$}
\psfrag{MA}[cc]{$M_A$} \psfrag{MB}[cc]{$M_B$}
\psfrag{MbA}[cc]{$M_{\bar{A}}$}
 \psfrag{MbarA}[cc]{$M_{\bar{A}}$}
\psfrag{MbarB}[cc]{$M_{\bar{B}}$}
 \psfrag{MbB}[cc]{$M_{\bar{B}}$}
\includegraphics[width=5in]{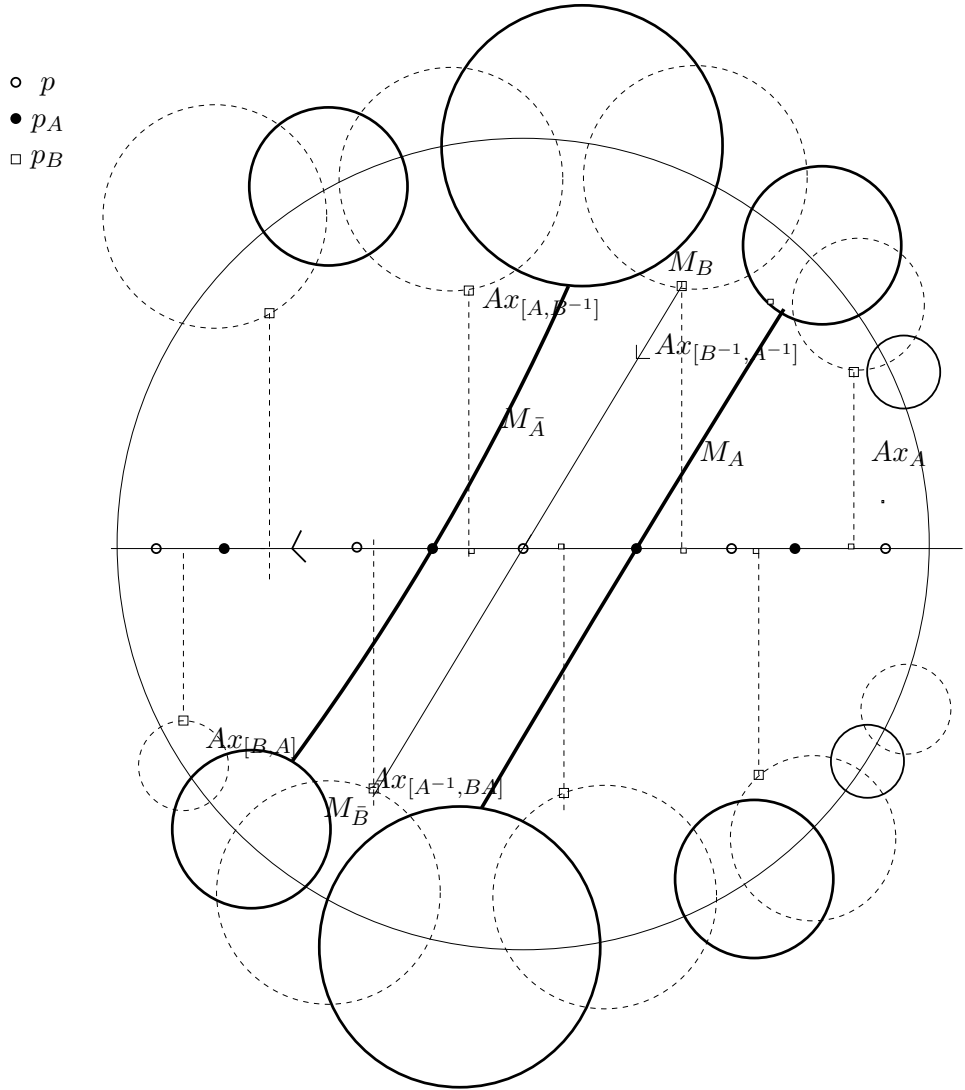}
\caption{The construction of the stopping generators: The
commutator axes are dark circles, $Ax_A$ and $Ax_B$ intersect at
the origin. The fundamental domain $F_K$ is bounded by the dark
lines $M_A,M_{\bar{A}}$ and the dotted circles $M_B,M_{\bar{B}}$.
\label{figure:intax3}}
\end{figure}

\begin{lemma}
\label{lemma:points on M} Let $M_A \in \D$ be the mutual
orthogonal to $Ax_{[B^{-1},A^{-1}]}$ and $Ax_{[A^{-1},B]}$ and
$M_B$  the mutual orthogonal to $Ax_{[B^{-1},A^{-1}]}$ and
$Ax_{[A,B^{-1}]}$. Then $M_A$ passes through $p_A$ and $M_B$
passes through $p_B$. Moreover, if $p_{\bar{A}}=A(p_A)$, and
$p_{\bar{B}}=B(p_B)$ and if $M_{\bar{A}}$ is the
 mutual orthogonal
  to $Ax_{[B,A]}$ and $Ax_{[A,B^{-1}]}$ and $M_{\bar{B}}$
is the mutual orthogonal to $Ax_{[B,A]}$ and $Ax_{[A^{-1},B]}$
then $M_{\bar{A}}$ passes through $p_{\bar{A}}$ and $M_{\bar{B}}$
passes through $p_{\bar{B}}$.
\end{lemma}

\begin{proof} (See figure~\ref{figure:intax3})
 We give the proof for $M_A$.  The others follow in
the same way. We factor $A$ and $B$ as
$A=E_pE_{p_{A}}=E_{p_{\bar{A}}}E_p$ and
$B=E_pE_{p_{B}}=E_{p_{\bar{B}}}E_p$. Then we easily compute
$$[B^{-1},A^{-1}]=E_{p_A}[A^{-1},B]E_{p_A}$$
so that $E_{p_A}$ interchanges the axis of $[B^{-1},A^{-1}]$ and
the axis of $[A^{-1},B]$.
 The elliptic sends
any geodesic lying on  $\D$ and passing through $p_A$ to its
inverse
 and moves any geodesic not passing through $p_A$.  Since
it interchanges the axes, it must
 send their
the mutual orthogonal to itself. Thus $p_A$ is on this orthogonal
as claimed. In the case of  $M_B$, note that $M_B$ is the same as
$\delta$.
\end{proof}

This construction yields a new description of a fundamental domain
for the intersecting acting case, one not previously used in the literature.

\begin{prop} \label{lemma:defM_A} The domain $F$
bounded by the geodesics $M_A$,$M_B$,$M_{\bar{A}}$ and
$M_{\bar{B}}$ is a fundamental domain for the action of $G$ on
$\D$ and $\D/G$ is a torus with a hole.  The domain $F_K$ obtained
by truncating $F$ along the commutator axes is a fundamental
domain for $G$ acting on the Nielsen  region.
\end{prop}

\begin{figure}[hbt]\centering
\psfrag{MA}[cc]{$\pi(M_A)$}
\psfrag{MB}[cc]{$\pi(M_B)$}
\psfrag{MbA}[cc]{$\pi(M_{\bar{A}})$}
\psfrag{MbB}[cc]{$\pi(M_{\bar{B}})$}
 \psfrag{PA}[cc]{$\pi(p_A)$}
\psfrag{PB}[cc]{$\pi(q_B)$}
 \psfrag{P}[cc]{$\pi(p)$}
\psfrag{PAxA}[cc]{$\pi(Ax_A)$}
 \psfrag{PAxB}[cc]{$\pi(Ax_B)$}
\includegraphics[width=5in]{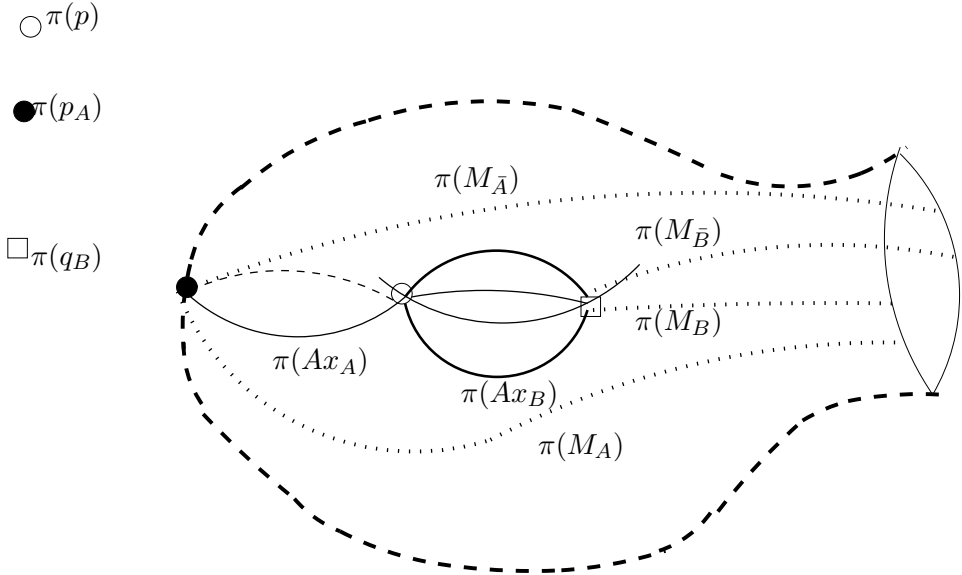}
\caption{Intersecting axes: the quotient $\D/G$
\label{figure:hole}}
\end{figure}

\begin{proof}
Note that the lines $M_A,M_B,M_{\bar{A}}$ and $M_{\bar{B}}$  are
mutually disjoint since if two of them intersected, they would
form, together with one of the commutator axes, a triangle with
two right angles.

Next, by definition, $M_A$ and $M_{\bar{A}}$ are identified by $A$
and $M_B$ and $M_{\bar{B}}$ are identified by $B$ and so, by
Poincar\'e's theorem, they bound a fundamental domain for $G$
acting on $\D$.

Let $D(I_{[B,A]})$ be the half plane with boundary $Ax_{[B,A]}$
that does not contain the point $p$. It is invariant under
$[B,A]$. Tiling $\D$ with copies of $F$, we see that the only
images of $F$ that intersect $D(I_{[B,A]})$ are of the form
$[B,A]^n(F)$, for some integer $n$.  It follows that
$D(I_{[B,A]})$ is stabilized by the cyclic group $\langle
[B,A]\rangle$ and the interval $I_{[B,A]}$ joining the fixed
points of $Ax_{[B,A]}$ is an interval of discontinuity for $G$.
The same is true for the other three intervals $I_{[A,B^{-1}]}$,
$I_{[A^{-1},B]}$ and $I_{[A^{-1},B^{-1}]}$ corresponding to the
other three commutators whose axes intersect $F$.

  Identifying the sides of $F$ we obtain a torus
with one hole. (See figure~\ref{figure:hole}). The ideal boundary
of the hole is made up the four segments of $\bar{F} \cap
I_{[B,A]}$, $\bar{F} \cap I_{[A,B^{-1}]}$, $\bar{F} \cap
I_{[A^{-1},B]}$ and $\bar{F} \cap I_{[A^{-1},B^{-1}]}$.
  There is one
funnel,  the quotient $D(I_{[B,A]})/\langle [B,A] \rangle$.

Truncating $F$ along the four commutator axes intersecting it we
have a fundamental domain for the Nielsen region of $G$ as
claimed.
\end{proof}

\subsection{Weierstrass points of the Schottky double}
\label{sec:wps schottky double int}

We note that the Schottky double $S$ is a surface of genus two
with one  funnel whose central curve is the projection of
$I_{[B,A]}$.  We now characterize its Weierstrass points.

   Let
 $p,p_A,p_B \in \D$ be the intersection points of
axes as in the previous section and let $q,q_A,q_B $ be their
respective reflections in $\partial\D$.
 Let $E$ be the elliptic of order
two with fixed point $p$.  $E$ leaves $\partial\D$ invariant and
its other fixed point is $q$. We can factor   $A$  as $A= E \circ
E_A$ where $E_A$ is the elliptic of order two with fixed point
$p_A$ and second fixed point $ q_A $.

Similarly factor $B=E \circ E_B$ where $E_B$ is a third elliptic
of order two with fixed points $p_B$ and  $ q_B $.

\begin{figure}[hbt]\centering
\psfrag{MA}[cc]{$\pi(M_A)$} \psfrag{MB}[cc]{$\pi(M_B)$}
\psfrag{MbA}[cc]{$\pi(M_{\bar{A}})$}
\psfrag{MbB}[cc]{$\pi(M_{\bar{B}})$} \psfrag{PA}[cc]{$\pi(p_A)$}
\psfrag{PB}[cc]{$\pi(p_B)$} \psfrag{P}[cc]{$\pi(p)$}
\psfrag{PAxA}[cc]{$\pi(Ax_A)$} \psfrag{PAxB}[cc]{$\pi(Ax_B)$}
\psfrag{Q}[cc]{$\pi(q)$} \psfrag{QA}[cc]{$\pi(q_A)$}
\psfrag{QB}[cc]{$\pi(q_B)$}
\includegraphics[width=5in]{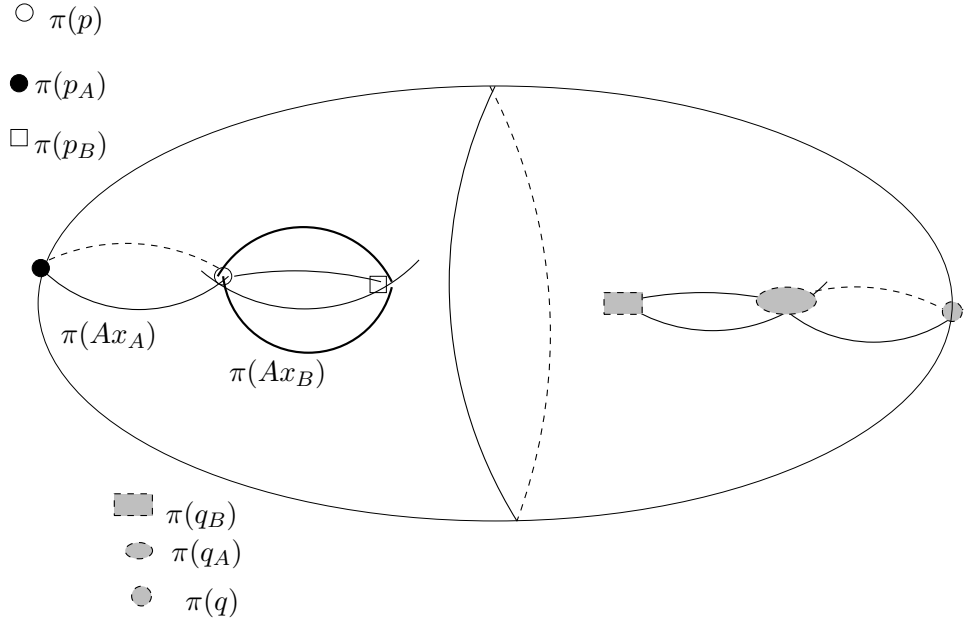}
\caption{Intersecting axes, one view of $\Omega(G)/G$: on the left
we
 have the projection of $\Di/G$ with the projections of the axes of
 $A$ and   $B$ and the Weierstrass points $\pi(p), \pi(p_A), \pi(p_B)$
indicated and on the right  we   have the projection of $\De/G$
with
 the
projections of the axes of $A$ and  $B$ and the Weierstrass points
$\pi(q), \pi(q_A), \pi(q_B)$ indicated. \label{figure:hole2}}
\end{figure}

\begin{figure}[hbt]\centering
\psfrag{MA}[cc]{$\pi(M_A)$}
 \psfrag{MB}[cc]{$\pi(M_B)$}
\psfrag{MbA}[cc]{$\pi(M_{\bar{A}})$}
\psfrag{MbB}[cc]{$\pi(M_{\bar{B}})$}
 \psfrag{PA}[cc]{$\pi(p_A)$}
\psfrag{PB}[cc]{$\pi(q_B)$}
 \psfrag{P}[cc]{$\pi(p)$}
\psfrag{PAxA}[cc]{$\pi{Ax_A}$}
 \psfrag{PAxB}[cc]{$\pi(Ax_B)$}
\psfrag{Q}[cc]{$\pi(q)$}
 \psfrag{QA}[cc]{$\pi(q_A)$}
\psfrag{QB}[cc]{$\pi(q_B)$}
\includegraphics[width=5in]{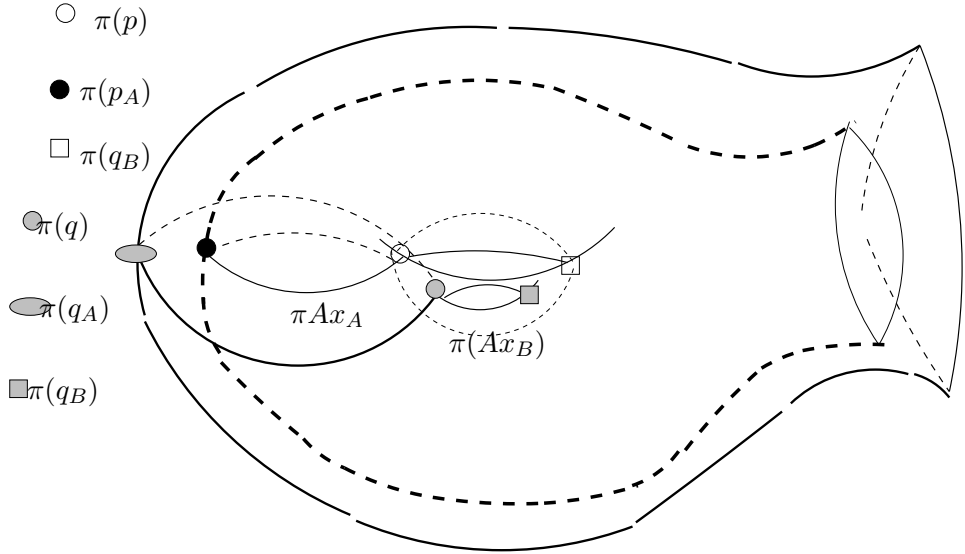}
\caption{Intersecting axes: alternate view of $\Omega(G)/G$. Here
  $\D/G$
is viewed as lying inside $J(\D/G)$ \label{figure:lost}}
\end{figure}
\begin{thm}
\label{thm:wps Schottky double int} Let $G=\langle A,B \rangle$ be
a discrete free Fuchsian group such that the axes of $A$ and $B$
intersect. Let $(p,q)$, $(p_{A},q_{A})$, $(p_{B},q_{B})$ be the
respective fixed points of the elliptics $E,E_A,E_B$ where
$A=E\circ E_A$ and $B=E \circ E_B$. Then the projections of these
points onto the surface $S$ of genus $2$ under the projection $
\Omega(G) \rightarrow \Omega(G)/G=S$,
  are the
Weierstrass points of $S$.
\end{thm}

\begin{proof} (Figures 8 and 9) We saw above that $\D/G$ is a torus with a hole so
that  its double $S$ is a surface of genus two with one funnel.

 Conjugation
of the group $G$ by $E$ sends $A$ to $A^{-1}$ and $B$ to $B^{-1}$
since $p$ lies on both their axes. It follows that $E$ induces  a
conformal involution $j$ of $S$. Moreover, since $E_A=E \circ A$
and $E_B=E \circ B$, these
 elliptics induce the same  involution on $S$. The projection of
 each of the six points is fixed by  $j$,  characterizing $j$
as the hyperelliptic involution.
\end{proof}

 There is again an anti-conformal involution $J$ of $S$
 induced by the symmetry of the
 group $G$ with respect to $\partial\Delta$, that is by $R_{\partial
 \Delta}$.  As a Schottky double, $S$ has a single funnel and the central geodesic
  of this
 funnel
is the projection, $\gamma$,  of
 $\Omega \cap \partial\Delta$; it is point-wise invariant under
 $J$.
 Note
 that $\gamma$  does not contain any of the Weierstrass points.

Both of the involutions $j$ and $J$ of $S$  fix  $\gamma$. The
first is conformal and fixes no point on $\gamma$. The other
involution $J$ is anti-conformal and fixes the entire curve point
wise.

There is a third involution $\hat{J}=Jj$.  This is the
anti-conformal map  on $S$ induced by the product $R_{\partial \D}
\circ E$.  Since $j$ fixes the Weierstrass points, both $J$ and
$\hat{J}$ have the same action on the Weierstrass points. Namely,
they both interchange the projections of $p$ and $q$, $p_A$ and
$q_A$ and $p_B$ and $q_B$ respectively.

\section{The convex hull and the convex core}
\label{sec:convexi}

 As in section \ref{sec:convexd}, we again
 consider the group $G$ acting on $\PP$, $\PPi$ and $\PPe$ which
 contain respectively the convex hull  $\CCC(G)$ in $\IH$ and its two
 boundaries $\CCC_{int}$ and $\CCC_{ext}$.

\subsection{Fundamental domains in $\PP$ and $\CCC$}
\label{sec:int axes in h3}  We now work in $\PP$ and draw in the
axes $Ax_A,Ax_B$, $Ax_{AB},Ax_{BA}$, $Ax_{A^{-1}B}$,
$Ax_{AB^{-1}}$; again we
 denote the intersection point of $Ax_A$ and $Ax_B$ by $p$, the
intersection point of $Ax_A$ and $Ax_{A^{-1}B}$ by $p_A$ and the
intersection point of $Ax_B$ and $Ax_{A^{-1}B}$ by $p_B$.

Denote the  half turns about lines orthogonal to $\PP$ that pass
through  $p,p_A,p_B$ by $H_p,H_{p_A},H_{p_B}$ respectively. When
restricted to $\PP$ these are the same as the rotations of order
two we  denoted by $E_p$, $E_{p_A}$ and $E_{p_B}$.

We  define the lines $M_A,M_{\bar{A}}$, $M_B,M_{\bar{B}}$ just as
we did in $\D$ and obtain a domain $F_{\PP}$ in $\PP$. The
arguments of the proof of lemma~\ref{lemma:defM_A} applied to
$\PP$ instead of $\D$ give an immediate proof of

\begin{prop}\label{prop:funddomain int} The domain
$F \subset \PP$ bounded by
the lines $M_A$,$M_{\bar{A}}$, $M_B$, $M_{\bar{B}}$ is a
fundamental domain for $G$ acting on $\PP$.  The truncation
$F_{\CCC}$ of $F$ along the axes of the commutators is a
fundamental domain for the convex hull $\CCC(G)$ in $\IH$.
\end{prop}

\begin{proof} Since $\CCC(G) \subset \PP$, it is precisely the
Nielsen convex region for $G$ acting on $\PP$ and
lemma~\ref{lemma:defM_A} applies. \end{proof} Again we define
$F_{int}=F_{\CCC} \cap \PPi$ and $F_{ext}=F_{\CCC} \cap \PPe$  as
domains in $\PPi$ and $\PPe$ respectively.  As an immediate
corollary we have

\begin{cor}
\label{cor:Fis corebdyint} The domain $F_{int} \cup F_{ext}$ is a
fundamental domain for $G$ acting on $\CCC_{int} \cup \CCC_{ext}
=\partial\CCC$.
\end{cor}

\section{The Pleated Surface}
\label{sec:plintaxes}

In this section we explicitly construct the pleated surface
$(\tilde{S},pl)$ for the intersecting axis case.  Again, we do
this by first constructing a Fuchsian group $\Gamma$ that
uniformizes the Nielsen double $S_K=\Delta/\Gamma$. We then
construct a pleating map $PL: \Delta \to \IH$ that intertwines the
actions of $\Gamma$ on $\Delta$ and $G$ on $\CCC \in \IH$.
Finally, we take quotients.

\subsection{The Fuchsian group for the Nielsen double}
\label{sec:group gamma int}  We begin with the  domain $F_K
\subset \D$ defined in proposition~\ref{lemma:defM_A} for the
Nielsen  region $K$.
 For readability, denote the reflection $R_{Ax_{[B^{-1},A^{-1}]}}$ by $R_0$.
 Define the
domain $F_{\Gamma}= F_K \cup R_0(F_K) \subset \Delta$. Let
$A'=R_0AR_0$ and $B'=R_0BR_0$ and let $\Gamma=\langle
A,B,A',B'\rangle$.

\begin{figure}[hbt]\centering
\psfrag{AXA-1B}[cc]{$Ax_{[A^{-1},B]}$}
\psfrag{AxB-1A-1}[cc]{$Ax_{[B^{-1},A^{-1}]}$}
\psfrag{AxAB-1}[cc]{$Ax_{[A,B^{-1}]}$}
 \psfrag{MA}[cc]{$M_A$}
\psfrag{MB}[cc]{$M_B$}
 \psfrag{MBA}[cc]{$M_{\bar{A}}$}
\psfrag{MB}[cc]{$M_{\bar{B}}$}
 \psfrag{M'bA}[cc]{$M'{\bar{A}}$}
\psfrag{M'bB}[cc]{$M'_{\bar{B}}$}
 \psfrag{A}[cc]{$A$}
\psfrag{B}[cc]{$B$}
 \psfrag{A'}[cc]{$A'$}
 \psfrag{B'}[cc]{$B'$}
\psfrag{MbA}[cc]{$M_{\bar{A}}$}
 \psfrag{MbarB}[cc]{$M_{\bar{B}}$}
\psfrag{MbarA}[cc]{$M_{\bar{A}}$}
 \psfrag{MbB}[cc]{$M_{\bar{B}}$}
\psfrag{R0(AxBA)}[cc]{$R_0(Ax_{[B,A]})$}
\psfrag{Ax[A,B^{-1}]}[cc]{$Ax_{[A,B^{-1}]}$}
\psfrag{AXBA}[cc]{$Ax_{[B,A]}$}
\psfrag{AXB-1A-1}[cc]{$Ax_{[B^{-1},A^{-1}]}$}
\psfrag{AXA-1B}[cc]{$Ax_{[A^{-1},B]}$}
\includegraphics[width=5in]{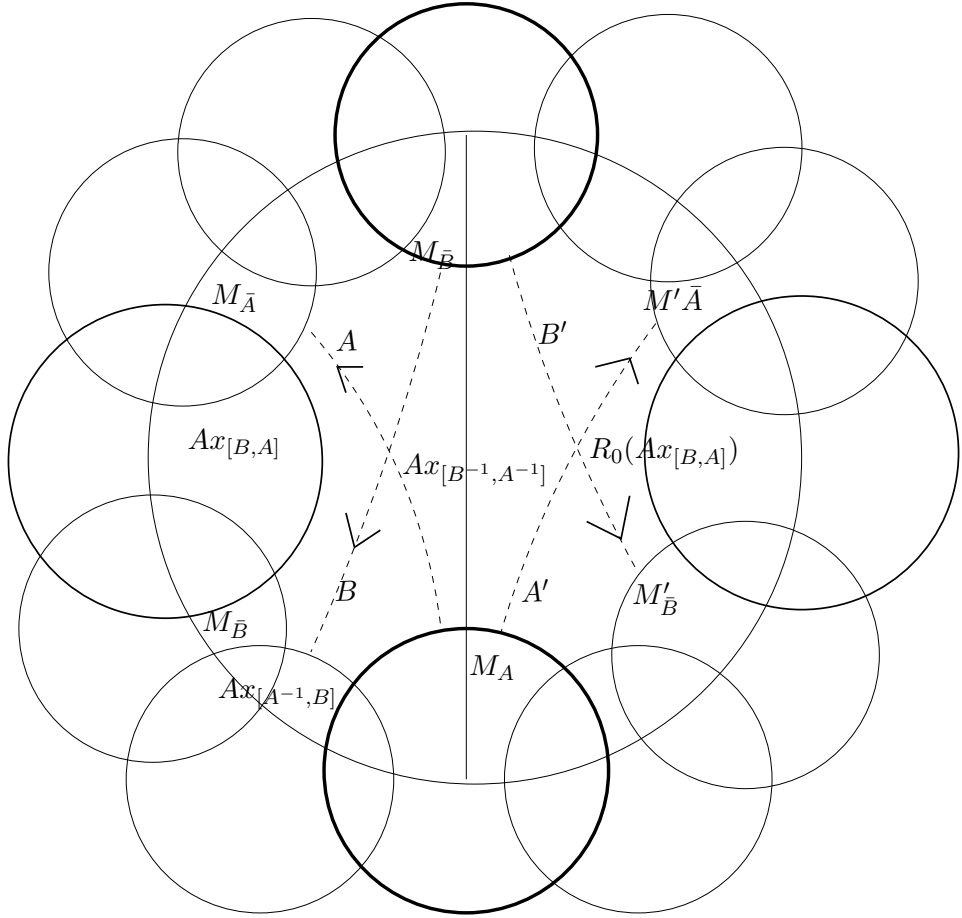}
\caption{The  domain $F_K$  and its reflection in
$Ax_{[B^{-1},A^{-1}]}$ form a fundamental domain for $\Gamma$}
\label{figure:reflint}
\end{figure}

\begin{prop} The domain $F_{\Gamma}$ is a fundamental domain for
$\Gamma$ acting on $\Delta$. The Fuchsian surface $S_K=\D/G$ is a
Riemann surface of genus two, the Nielsen double.
\end{prop}

\begin{proof} (See figure~\ref{figure:reflint}) The sides of the domain $F_K$ of
proposition~\ref{lemma:defM_A} are $M_A$,$M_{\bar{A}}$,$M_{B}$ and
$M_{\bar{B}}$.
 Set $M_A'=R_0(M_A)$, $M_B'=R_0(M_B)$,
$M_{\bar{A}}'=R_0(M_{\bar{A}})$ and
$M_{\bar{B}}'=R_0(M_{\bar{B}})$. The axis $Ax_{[B^{-1},A^{-1}]}$
divides $M_A=M_A'$ into two segments; without causing confusion we
will say $M_A$ is the segment in $F_K$ and $M_A'$ is the segment
in $R_0(F_K) $ and similarly for $M_B$ and $M_B'$.

Using the side identifications
$$A:M_A \rightarrow M_{\bar{A}}, \, B:M_B \rightarrow
M_{\bar{B}}$$ we see that if we set $A'= R_0 A R_0$ and $B'=R_0 B
R_0$ we obtain the identifications
$$A':M_A' \rightarrow M_{\bar{A}}', \,  B':M_B' \rightarrow
M_{\bar{B}}'$$

Now denote  the reflection in the axis $Ax_{[A^{-1},B]}$ by $R_1$
and note that $R_1=BR_0B^{-1}$.  The axes $Ax_{[A^{-1},B]}$ and
$Ax_{[B^{-1},A^{-1}}$ are both orthogonal to $M_A$. Thus we see
that
$$R_1R_0=(B')^{-1}B:Ax_{[A,B^{-1}]}\rightarrow
R_0(Ax_{[A,B^{-1}]})$$ Similarly
$$(A')^{-1}A:Ax_{[A^{-1},B]}\rightarrow
R_0(Ax_{[A^{-1},B]})$$
$$A'B'B^{-1}A^{-1}:Ax_{[B,A]}
\rightarrow R_0(Ax_{[B,A]})$$

Including the two vertices where $M_A$ meets $M_A'$ and where
$M_B$ meets $M_B'$ there are $14$ vertices in all.  They fall into
two cycles of four vertices where the angles are all right
 and two cycles of three
vertices, two where there are right angles and one where there is
a straight angle.
That $\Gamma$ is Fuchsian and
$F_{\Gamma}$ is a fundamental domain now follows from Poincar\'e's
theorem. Identifying  the sides $M_A$ and $M_{\bar{A}}$ and the
sides $M_B$ and $M_{\bar{B}}$
 of $F_{K}$ yields a
torus with a hole where the remaining unidentified sides (arcs of
commutator axes) fit together to form the boundary of the hole.
Similarly, identifying corresponding sides of $R_0(F_{K})$ yields
another torus with a hole.  In $F_{\Gamma}$, the commutator axis
sides are identified and the two tori with a hole join up to form
the Nielsen double $S_K$ as a compact surface   of genus two.
\end{proof}

\subsection{The pleating map}
\label{sec:pleating map int}
 We now construct the  pleating map
$pl: \D/\Gamma \rightarrow \IH/G$ just as we did in
section~\ref{sec:pleating map nonint}. We use the definitions of
$F_{\Gamma} \subset \D$ as a union of two copies of $F_K$, and the
the domains  $F_{int} \subset \PP_{int}$ and $F_{ext} \subset
\PP_{ext}$ from section~\ref{sec:int axes in h3}. Again we begin
by defining a pleating map $PL:\D \rightarrow \IH$. We first
define $PL$  on $F_{\Gamma}$ and then extend by the groups
$\Gamma$ and $G$.

Set  $PL: F_K \rightarrow F_{ext}$. Next we set
 $PL: R_{0}(F_K) \rightarrow F_{int}$.
 Define the group homomorphism $\phi: \Gamma \rightarrow G$ by first
defining  it on the generators:
$$A \mapsto A, B \mapsto B, A' \mapsto  id, B' \mapsto id.$$
 To show that this map on generators gives a group homomorphism
from $\Gamma$ to $G$ we must show that $\phi$ preserves the
defining relation(s)
 of $\Gamma$. As in the previous case, we know that there are no
  reflection relations and the cycle relations
form  a complete set of relations for $\Gamma$. Again, it is easy
to calculate that substituting $\phi(A)$, $\phi(A')$, $\phi(B)$,
 and $ \phi(B')$ into the cycle relation gives the identity.

We now extend $PL$ to a map from $\D$ into $\IH$. That is, for $x
\in \D$ let $g \in \Gamma$ be chosen so that $g(x) \in F_{\Gamma}$
and set $PL(x) = \phi(g^{-1}) \circ PL \circ g(x)$. This is well
defined  when  $g(x)$ lies interior to $F_{\Gamma}$ and when
$g(x)$ lies on a boundary curve of $F_{\Gamma}$ or on $Ax_{[B,A]}$
it is defined by continuity. Thus,  $(\D,PL)$ is a pleated surface
with image in $\IH$ whose pleating locus is the image of the axis
sides of $F_{\Gamma}$ and all their images under $\Gamma$.

We now define the map $pl:S_K \rightarrow \IH/G$ by taking
quotients.   The map $pl$ is clearly a hyperbolic isometry since
$PL$ is and its pleating locus consists of the images of  the
commutator axes.

By corollary~\ref{cor:Fis corebdyint}, we can apply the group $G$
to identify $PL(F_{\Gamma})$ with $\partial\CCC$, and taking
quotients, identify $pl(S_K)$ with the convex core boundary
$\partial{N}=\CCC_{int}/G\cup \CCC_{ext}/G$.

\subsection{The Weierstrass points of the Nielsen double $S_K$}
\label{sec:wps nielsen double int}

We again now determine the Weierstrass points of $S_K$ in a manner
reminiscent of our  discussion of the Schottky double. The map
$E_p$ is a conformal involution that maps $F_K$ to itself and maps
$R_0(F_K)$ to its image under $AB(B')^{-1}(A')^{-1}$. Therefore
$E_p$ preserves the $F_{\Gamma}$ tiling of $\D$ which implies that
$E_p$ conjugates $\Gamma$ to itself.
 It follows
that $E_p$ induces a conformal involution $j$ of $S_K$. It is easy
to check that the points $p,p_A,p_B$, $p',p_A',p_B'$ are mapped by
$E_p$ to  points that are equivalent under the action of $\Gamma$
and thus project to  fixed points of $j$. This characterizes $j$
as the hyperelliptic involution of $S_K$.

Note that $PL(p) \in F_{int}$ and $PL(p') \in F_{ext}$, but as
points in $\PP$, $PL(p)=PL(p')$ and similarly for the other pairs
of Weierstrass points.

For the other elliptics fixing these six points in $F_{\Gamma}$ we
have:
$$E_{p_A}=E_pA, \, E_{p_B}=E_pB, \,
E_{p'}=E_pAB(B')^{-1}(A')^{-1},$$ $$E_{p_A'}=E_{p'}A', \,
E_{p_B'}=E_{p'}B'$$ so that they also induce the involution $j$.

\section{Weierstrass points and lines in the handlebody}

\label{sec:gen weierstrass int} In this section we give an
explicit description of the handlebody $H$ as $\cal{S} \times I$
in the case where the axes of the generators intersect.  We show
that again, $H$ admits a unique order two isometry $j$ that we
call the {\em hyperelliptic isometry of $H$}. It fixes six unique
geodesic line segments, which we call the {\sl Weierstrass lines}
of $H$.

 We show, again, that $H$ is
foliated by surfaces $\cal{S}(s)=(\cal{S},s)$ that are at distance
$s$ from the convex core, and that each $\cal{S}(s)$ is fixed by
$j$. The restriction of $j$ to  $\cal{S}(s)$ fixes six {\sl
generalized Weierstrass points} on $\cal{S}(s)$. Moreover, on each
$\cal{S}(s)$  we
 find an orientation reversing involution $J$ that
interchanges pairs of   generalized Weierstrass points.

\subsection{Construction of the foliation}
\label{sec:equisurface int}

We modify our construction in  section~\ref{sec:gen weierstrass
non int}  to construct the {\sl hyperelliptic isometry for $H$} as
follows.  In section~\ref{sec:int axes in h3}, we defined the
half-turns about lines orthogonal to $\PP$
at the points $p,p_A,p_B$ in $\PP$.  Set $p_{int}$ and $p_{ext}$
as the point $p$ considered in $\PPi$ or $\PPe$ respectively, and
similarly for the other points.

Following our notation in section~\ref{sec:equisurface disjoint}
we parameterize these orthogonal lines by arc length and denote
them by $V_p(t)$,$V_{p_A}(t)$ and $V_{p_B}(t)$ respectively, such
that $V_p(0)=p,$ $V_{p_A}(0)=p_A$ and $V_{p_B}(0)=p_B$. Direct
them so that for $t>0$, $V_p(t)$ is in the half-space with
boundary $\PPe$ and similarly for the other lines.  With this
convention we write $V_p(\pm s)$,$V_{p_A}(\pm s)$ and $V_{p_B}(\pm
s)$ where $s=|t|$.

Since the fundamental domains $F_{ext}$ and $F_K$ are related by
orthogonal projection we see that for $p,q \in F_K$ we have
$p_{ext}=PL(q) \in F_{ext}$ is equal to $\lim_{s \to
0}V_p(s)=V_p(0^+)$ and $p_{int}=PL(p) \in F_{int}$ is equal to
$\lim_{s \to 0}V_p(-s)=V_p(0^-)$. Also, $p \in F_K \subset \Di$ is
equal to $\lim_{s \to \infty}V_p(-s)$ and $q \in
R_{\partial\D}(F_K) \subset \De$ is equal to $\lim_{s \to
\infty}V_p(s)$ and similarly for the other points.

 We saw that $A$ and $B$ may be factored into
products of pairs of these half turns and it follows that $G$ is a
normal subgroup of index two in the group  generated by these half
turns. Since the half-turns all differ by
elements of $G$, their actions by conjugation on $G$ all induce
the same order two automorphism $j$ under the projection $\pi:\IH
\rightarrow \IH/G=H$. The fixed points of $j$ in $H$ are precisely
the points on the images of the six lines $V_p(\pm
s)$,$V_{p_A}(\pm s)$ and $V_{p_B}(\pm s)$: denote them by
$\pi(V_p)(\pm s)$,$\pi(V_{p_A})(\pm s)$ and $\pi(V_{p_B})(\pm s)$.

We work with the lines $V_p$ and $\pi(V_p)$ but we have analogous
statements for the  lines $V_{p_A}$ and $V_{p_B}$ and their
projections. As we saw in theorem~\ref{thm:wps Schottky double
int}, the endpoints $p \in \Di$ and $q \in \De$ of $V_p$ project
to Weierstrass points on the Schottky double
$S=\partial{H}=\cal{S}(\infty)$.

 In theorem~\ref{thm:wps of Nielsen double nonint} we
found the Weierstrass points of the Nielsen double $S_K$. In the
fundamental domain $F_{\Gamma} \subset \D$ for $\Gamma$ we found
points $p,q \in F_{\Gamma} \in \D$ such that
$p_{int}=PL(p)=V_p(0^-)$ and $p_{ext}=PL(q)=V_p(0^+)$.  For the
moment, call the hyperelliptic involution on $S_K$, $j_{S_K}$ and
the hyperelliptic involution on $H$, $j_H$. By construction, the
pleating map $pl$ commutes with these involutions: $pl \circ
j_{S_K}=j_H \circ pl$. It follows that the projections $\pi(PL(p))
\in N_{int}$ and $\pi(PL(q)) \in N_{ext}$ are fixed under $j_H$.
By the same argument, we find two other  pairs of points
$\pi(PL(p_A)),\pi(PL(q_A))$ on $\pi(L_A)$ and
$\pi(PL(p_B)),\pi(PL(q_B))$ on $\pi(L_B)$ fixed by $j_H$. They lie
on the genus two surface $\partial{N}=\cal{S}(0)$ and are its {\em
generalized Weierstrass points}.

Each of the segments of the line $\pi(V_p(\pm s)$  is point-wise
fixed under $j_H$. We call each a {\em generalized Weierstrass
line of H}. Similarly, we have two other pairs of generalized
Weierstrass lines $\pi(V_{p_A}(\pm s))$, and $\pi(V_{p_B}(\pm
s))$.

  We next
 construct the family of surfaces $\cal{S}(s)$. We will show that
 $j_H$, which we denote again simply by $j$ since it will not
 cause confusion, is an order two isometry of $\cal{S}(s)$ that has as its set of fixed points,
  the  intersection points of these six lines with the surface.

We again consider the family of equidistant surfaces  $\Pi(\pm s)$
with boundary $\partial\Delta$ at distance $s$ from the plane
$\PP$. Using orthogonal projection, we can project ${\CCC_{ext}}$
onto $\Pi(s)$ to obtain $\CCC_{ext}(s)$ and project
$\CCC_{int}$ onto $\Pi(-s)$  to obtain $\CCC_{int}(s)$.

Each quotient $N_{ext}(s)= \CCC_{ext}(s) /G$ and $N_{int}(s) =
\CCC_{int}(s)/G$ is topologically a torus with one boundary curve.
In analogy with our construction of the Schottky double, we want
to join these boundary curves  by  a funnel to form the surface
$\cal{S}(s)$.

Let  $Eq_{[B,A]}(s)$ be the equidistant cylinder in $\IH$ about
$Ax_{[B,A]}$. Note that, for each $s>0$,  $[B,A]$ maps
$Eq_{[B,A]}(s)$ to itself. Because $Ax_{[B,A]}$ lies in $\PP$,
$[B,A](\pm s)= Eq_{[B,A]}(s)\cap \Pi(\pm s)$ is the orthogonal
projection of $Ax_{[B,A]}$ to $\Pi(\pm s)$. The curves $[B,A](\pm
s)$ intersect $\CCC_{ext}(s)$ and $\CCC_{int}(s)$ in their
boundary curves. The curves $[B,A](\pm s)$ divide $Eq_{[B,A]}(s)$
 into
two pieces, each an infinite topological strip. Both of these
strips intersect $\PP$, one on each side of the axis of $[B,A]$.
Thus one strip intersects $\CCC$ and one does not; we denote the
strip disjoint from $\CCC$ by $Rec_{[B,A]}(s)$.

Similarly, we form strips for each of the boundary curves of
$\CCC$. The quotient of each of these strips under the action of
$G$ is a funnel; there is only one conjugacy class of boundary
curves, so only one funnel, denoted $F_{[B,A]}(s)$.

The  surface $\cal{S}(s)$ is defined as
$$\cal{S}(s)=N_{ext}(s) \cup N_{int}(s) \cup F_{[B,A]}(s) $$
 We have
\begin{prop} \label{prop:calS is invariant int} The surface
$\cal{S}(s)$ is invariant under the involution $j$. \end{prop}
\begin{proof} By construction, the distance from each point of $\cal{S}(s)$ to
the convex core $N$ is $s$.  As the involution $j:H \rightarrow H$
is an isometry and it preserves $N$, it leaves $\cal{S}(s)$
invariant. \end{proof}

We want to find the fixed points of $j$ acting on $\cal{S}(s)$. To
this end consider the intersection points of the line segments
$V_p(\pm s)$  with $\Pi(\pm s)$;  denote them respectively by
$q(s)$ and $p(s)$.  Similarly, define the points $q_A(s)$ and
$p_A(s)$ as the intersection points of the line segments
$V_{p_A}(\pm s)$ with $\Pi(\pm s)$ and define the points $q_B(s)$
and $p_B(s)$ as the intersection points of the line segments
$V_{p_B}(\pm s)$  with $\Pi(\pm s)$.

\begin{prop} The projections of the points $p(s),q(s)$,
$p_A(s),q_A(s)$ and  $p_B(s),q_B(s)$ lie on the surface
$\cal{S}(s) \subset H$ and comprise the set of fixed points of $j$
restricted to $\cal{S}(s)$.
\end{prop}
\begin{proof} Each of these points lies on one of the six
generalized Weierstrass lines in $H$ at distance $s$ from the
convex core. These lines are distinct so the points are distinct.
Since these lines are fixed point-wise by $j$, and since $j$
leaves $\cal{S}(s)$ invariant, these are fixed points of $j$ on
$\cal{S}(s)$ as claimed.
\end{proof}

We call these projected points the {\em generalized Weierstrass
points} of $\cal{S}(s)$.

We summarize these results as

\begin{thm}
\label{thm:foliation int} The projections of the
 line segments $V_p(\pm s)$, $V_{p_A}(\pm s)$ and $V_{p_B}(\pm s)$, orthogonal to
  $ \PP \subset H$ are point-wise fixed by the involution
 $j$.   The endpoints of each segment, lying in the boundary of
the convex core,  are the generalized Weierstrass points of
$(\cal{S},0)$. Each line is parameterized by hyperbolic arc length
and is called a generalized Weierstrass line. The points at
distance $s$ on the generalized Weierstrass lines are the
generalized Weierstrass points of the surface $(\cal{S},s)$. The
endpoints of the Weierstrass lines lie on the boundary the
handlebody and are the Weierstrass points of the boundary surface
$S=(\cal{S},\infty)$.
\end{thm}

\begin{figure}[hbt]\centering
\psfrag{R}[cc]{$r_{0}$} \psfrag{RR}[cc]{$\;\;r_{-}$}
\psfrag{RRR}[cc]{$\;\;r_{+}$} \psfrag{K}[cc]{$q_+$}
\psfrag{L}[cc]{$q_0$} \psfrag{M}[cc]{$q_{-}$}
\psfrag{N}[cc]{$p_{-}$} \psfrag{X}[cc]{$p_{0}$}
\psfrag{Y}[cc]{$p_{+}$}
\includegraphics[width=5in]{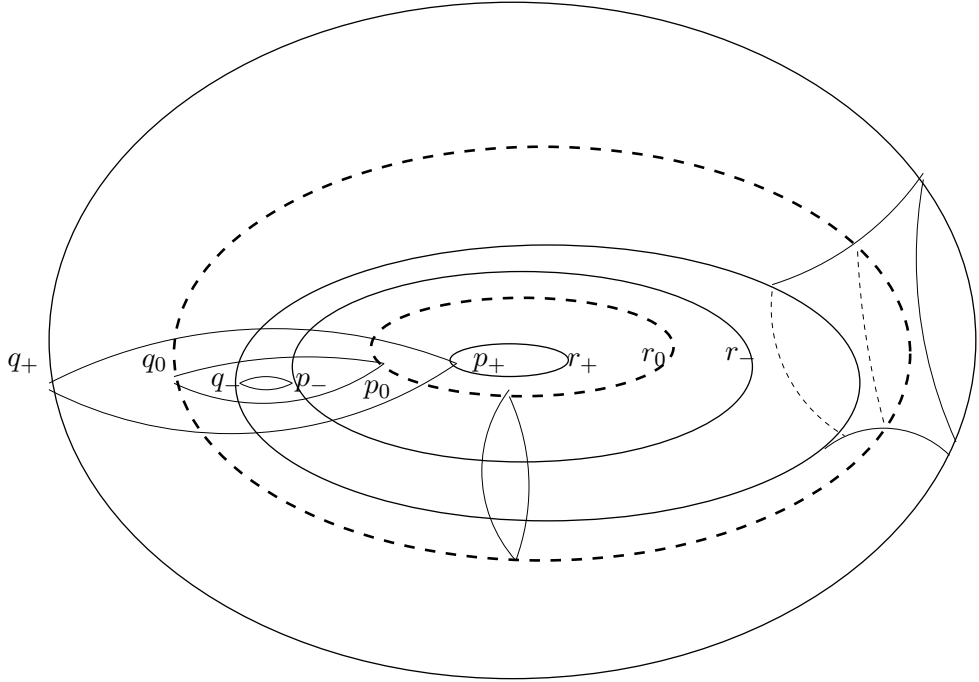}
\caption{The convex core is a one holed torus inside the
  handlebody. Its boundary is indicated in dotted lines. The
  Weierstrass points are marked as follows:
$p_{-}=\pi(V_p(-\infty))$, $p_0=\pi(V_p(0))$,
$p_{+}=\pi(V_p(+\infty))$, $q_{-}=\pi(V_{p_A}(-\infty))$,
$q_0=\pi(V_{p_A}(0))$,
  $q_{+}=\pi(V_{p_A}(+\infty))$,$r_{-}=\pi(V_{p_B}(-\infty))$;
  $r_0=\pi(V_{p_B}(0))$,
$r_{+}=\pi(V_{p_B}(+\infty))$. The Weierstrass lines join $q_{+}$
to $q_0$, $q_{0}$ to $q_{-}$, $p_{-}$ to $p_0$, $p_0$ to $p_{+}$,
$r_{-}$ to $r_0$, and $r_0$ to $r_{+}$. \label{figure:hole7}}
\end{figure}

Note that there are no generalized Weierstrass points on the
central curve of the funnel.  It is invariant under $j$ and is
mapped to its inverse.
\subsection{Anticonformal involutions}
\label{sec:anticonf int}
 We can also construct the anti-conformal
involution $J$ acting on $H$, by defining it to be the self-map of
$H$ induced by the anti-conformal reflection in the plane $\PP$,
$R_{\PP}$. Since the axes of the generators $A$ and $B$ of $G$ lie
in $\PP$, $R_{\PP}$ fixes theses axes point-wise and since
$R_{\PP}$ is an orientation reversing map, it induces an
orientation reversing map on $H$.

We want to see how $J$ acts on
$$\cal{S}(s)=N_{ext}(s) \cup N_{int}(s) \cup F_{[B,A]}(s) $$
Since $R_{\PP}$ is an isometry,
$R_{\PP}(\CCC_{ext}(s))=\CCC_{int}(s)$ and $J(N_{ext})=N_{int}$.
As $J$ is orientation reversing, it maps the projection of the
line $V_p(s))$ to the projection of the line $V_p(-s)$ and
interchanges the projections of $p(s)$ and $q(s)$.  It acts on the
other lines
 and generalized Weierstrass points similarly.

 The involution maps $Rec_{[B,A]}(s)$ to itself, fixing point-wise the curve that
is the intersection of the plane $\PP$ with $Rec_{[B,A]}(s)$. We
call the projection of this curve the central curve of the funnel.
Since $R_{\PP}$ interchanges the part of $Rec_{[B,A]}(s)$ on one
side of $\PP$ with the part on the other, we deduce that $J$ maps
the  funnel $F_{[B,A]}(s)$ to itself, interchanging its boundary
curves and sends the central curve to its inverse.

There is another  anti-conformal self-map  of $H$ induced by
$R_{\PP} \circ H_L$; $\hat{J}=J \circ j$.  To see how $\hat{J}$
acts on $\cal{S}(s)$ we first look at the fundamental domains
$F_{ext}(s)$ and $F_{int}(s)$:   $J$ interchanges these domains
but $j$ leaves them invariant, so $\hat{J}$ interchanges them.

  It follows that $\hat{J}$ maps $N_{ext}(s)$ to $N_{int}(s)$ maps
  the funnel to itself, interchanging the boundary curves
   and interchanges the generalized Weierstrass points.

\end{article}

\begin{thebibliography}{99}

\bibitem{Ahl}  Ahlfors, L.
\newblock {\sl Complex Analysis}
\newblock McGraw-Hill, third edition,  (1979).

\bibitem{Bers}  Bers, L.  {\sl Nielsen extensions of
Riemann surfaces}, Ann. Acad. Sci. Fenn. {\bf 2} (1976) 29-34.
%
%


\bibitem{Fench} Fenchel,  W. {\em Elementary Geometry in Hyperbolic Space},
de Gruyter Studies in Mathematics, $11$, Berlin-New York,(1989).


\bibitem{Gint}  Gilman, J. {\sl Inequalities and  Discreteness}
Canad. J. Math, {\bf  XL}(1)  (1988) 115-130.
\bibitem{Galg} Gilman, J.  {\sl Two generator discrete subgroups of $PSL(2,\RR)$},
Memoirs of the AMS, Vol. 117, No. 561 (1995).
\bibitem{GilKeen}  Gilman, J.  and  Keen, L. {\sl Word sequences and intersection numbers},
 Complex manifolds and hyperbolic
 geometry (Guanajuato, 2001), 231--249, Contemp. Math., 311, Amer. Math. Soc., Providence, RI, 2002.

\bibitem{GilMas}  Gilman, J. and  Maskit, B.  {\em An Algorithm for two-generator discrete
groups}, Mich. Math. J., {\bf 38 (1)} (1991), 13-32.


\bibitem{KACT}  Keen, L. {\em Canonical Polygons for Finitely Generated Fuchsian Groups.}  Acta
    Mathematica {\bf 115} (1966)

\bibitem{KANN} Keen,  L. {\em Intrinsic Moduli on Riemann Surfaces}.  Annals of Mathematics
{\bf 84} \#3 (1966),  404--420
\bibitem{Kap}  Kapovich, M.  {\sl Hyperbolic Manifolds and Discrete Groups}
Birkauser (2001).

\bibitem{Knapp}  Knapp, A. W.  {\sl Doubly generated Fuchsian groups}
Mich. Math. J. {\bf 15} (1968) 289-304.

\bibitem{Hur}  Hurwitz, A. {\em Algebraische Gebilde mit eindeutigen
Tranformationen in sic}, Math. Ann {\bf 41}, (1893) 409-442.
%
\bibitem{Mas} Maskit, B.  {\em Kleinian groups}, Springer-Verlag,
(1988).

\bibitem{Mat}  Matelski, P. {\em The classification of discrete two-generator
subgroups of $PSL(2, \RR)$}, Israel. J. Math. {\bf 42} (1982)
309-317.




\bibitem{PR} Purzitsky, N. and   Rosenberger, G.  {\em All two-generator Fuchsian
groups} Math. Z. {\bf 128} (1972) 245-251,  Correction: Math. Z. {\bf 132}
(1973) 261-262.

\bibitem{Pu}  Purzitsky, N. {\em All two generator Fuchsian groups} Math.Z. {\bf 147}(1976) 87-92.
\bibitem{RR}  Rosenberger, N.  {\em All generating pairs of all two-generator
Fuchsian groups}, Arch. Mat {\bf 46} (1986) 198-204.




\bibitem{Spr} Springer, G. {\em  Introduction to Riemann Surfaces},
Addison-Wessley,  (1957).

\end{thebibliography}
\end{document}